\documentclass{article}

\usepackage[english]{babel}

\usepackage[letterpaper,top=2cm,bottom=2cm,left=3cm,right=3cm,marginparwidth=1.75cm]{geometry}

\usepackage{amssymb}
\usepackage{amsthm}
\usepackage{amsmath}
\usepackage{empheq}
\usepackage{amscd}
\usepackage{graphicx}
\usepackage{graphics}
\usepackage[noadjust]{cite}
\usepackage{caption}
\usepackage{multirow}
\usepackage{array}
\usepackage{rotating}
\usepackage[norelsize,ruled]{algorithm2e}

\usepackage{indentfirst}
\usepackage{tikz}
\usepackage{calc}
\usepackage{epsfig}
\usepackage[numbers]{natbib}
\usepackage[makeroom]{cancel}
\usepackage{multicol}
\usepackage{csquotes}
\usepackage{enumitem}

\usepackage{hyperref}       
\usepackage{optidef}
\usepackage{nicefrac}

\newcommand{\oneD}{\text{\scriptsize{1D}}}
\newcommand{\opt}{\text{\scriptsize{opt}}}

\newcommand{\bD}{\bold{D}}
\newcommand{\be}{\bold{e}}
\newcommand{\bg}{\bold{g}}
\newcommand{\blam}{\boldsymbol{\lambda}}
\newcommand{\bm}{\bold{m}}
\newcommand{\bv}{\bold{v}}
\newcommand{\nrmg}{\text{nrm}(\bold{g})}
\newcommand{\nrm}{\text{nrm}}

\usepackage{epstopdf} 
\usepackage{mathptmx}      

\usepackage[title]{appendix}

\newtheorem{proposition}{Proposition}[section]

\title{Automatic balancing parameter selection\\for Tikhonov-TV regularization}
\author{Ali Gholami and Silvia Gazzola}

\begin{document}
\maketitle

\begin{abstract}
This paper considers large-scale linear ill-posed inverse problems whose solutions 
can be represented as sums of smooth and piecewise constant components. 
To solve such problems we consider regularizers consisting of two terms that must be balanced. Namely, a Tikhonov term guarantees the smoothness of the smooth solution component, while a total-variation (TV) regularizer promotes blockiness of the non-smooth solution component. A scalar parameter allows to balance between these two terms and, hence, to appropriately separate and regularize the smooth and non-smooth components of the solution. 
This paper proposes an efficient algorithm to solve this regularization problem by the alternating direction method of multipliers (ADMM). Furthermore, a novel algorithm for automatic choice of the balancing parameter is introduced, using robust statistics. The proposed approach is supported by some theoretical analysis, and numerical experiments concerned with different inverse problems are presented to validate the choice of the balancing parameter.
\end{abstract}

\section{Introduction}
This paper is concerned with the solution of ill-conditioned linear systems of equations of the form
\begin{equation}\label{main_eq}
\bold{d=Gm + e}, 
\end{equation}
arising from the discretization of continuous inverse problems in different fields of applied sciences and engineering; see, e.g., \cite{Aster_2004_PEI, Tarantola_2005_IPT}. Here $\bold{G}$ is an $M\times N$ real matrix representing a forward operator, $\bold{d}$ and $\bold{e}$ are real vectors of lengths $M$ representing known measured data and unknown Gaussian white noise, respectively, and the unknown $\bold{m}$ is a real vector of length $N$ representing a quantity of interest. 
The singular values of $\bold{G}$ may concentrate close to the origin and decay to zero, making the system ill-conditioned. 
Both one-dimensional (1D) and two-dimensional (2D) problems are considered in this paper; in the 2D case, the model parameters $\bold{m}$ are obtained by stacking the columns of a rectangular $N_z\times N_x$ array, so that $N=N_z\cdot N_x$. 
Due to the presence of noise and the ill-conditioning of the system, a naive (e.g., least squares) solution of \eqref{main_eq} possibly gives a meaningless estimate of the unknown 
$\bold{m}$. Thus we seek to determine an approximation of $\bold{m}$ by regularizing problem \eqref{main_eq}, i.e., by incorporating some a priori information about $\bold{m}$ into the problem formulation; see
\cite{Hansen_2010_DIP}. 

In this paper, we assume that the solution to the problem \eqref{main_eq} is piecewise smooth, e.g., it consists of a smooth background while containing regions with discontinuities. This allows us to decompose $\bold{m}$ as the sum of two components $\bold{m}_1$ and $\bold{m}_2$, i.e., 
\begin{equation}
\bold{m} = \bold{m}_1+\bold{m}_2,
\end{equation}
where $\bold{m}_1$ is the piecewise constant component of  $\bold{m}$ and $\bold{m}_2$ is the smooth component of  $\bold{m}$. Following \cite{Gholami_2013_BCT}, we incorporate regularizarization by considering a Tikhonov regularization term for the smooth component and a(n anisotropic) total variation (TV) term for the piecewise constant component.
A scalar balancing parameter $\beta>0$ is included to control the amount of regularization associated to each term and, hence, to properly separate the two complementary components of the solution.
Assuming that the value $\varepsilon=\|\bold{e}\|_2^2$ is available, the pair ($\bold{m}_1,\bold{m}_2$) is obtained by solving the following Tikhonov-TV constrained optimization problem
\begin{mini} 
{\bold{m}_1,\bold{m}_2}{\|\bold{D}_1\bold{m}_1\|_1+\frac{\beta}{2}\|\bold{D}_2 \bold{m}_2\|_2^2}
{\label{main_con}}{}
\addConstraint {\|\bold{G}(\bold{m}_1+\bold{m}_2)-\bold{d}\|_2^2}{=\varepsilon},
\end{mini}
where $\|\cdot\|_p$, $p=1,2$, denotes the vector $p$-norm, and 
$\bold{D}_1$ and $\bold{D}_2$ are scaled finite difference discretizations of the gradient and partial second derivative operators, 
respectively. 


The Tikhonov-TV regularization method \eqref{main_con} has gained considerable attention in scientific applications, for example: image denoising, inpainting and deblurring \cite{Gholami_2013_BCT}, signal decomposition \cite{huska2019convex}, optical flow computation \cite{kalmoun2020new}, signal filtering \cite{selesnick2014simultaneous}, seismic travel time tomography \cite{alrajawi2017inversion}, and seismic full waveform inversion \cite{Aghamiry_2019_CROb}.
This paper proposes an alternative algorithm to the one in \cite{Gholami_2013_BCT}, featuring a fully automatic balancing parameter choice strategy, which can be reliably used to solve large scale problems arising in the applications mentioned above. Note that, in \cite{Gholami_2013_BCT}, $\xi=\beta^{-1}$ was used as balancing parameter; using $\beta$ is more convenient for the approach presented in this paper. 

In this paper we solve \eqref{main_con} using the alternating direction method of multipliers (ADMM) \cite{Gabay_1976_ADA,chen1994proximal}.  This method has proved to be particularly advantageous for dealing with large scale problems such as those arising in geophysics and other scientific applications \cite{Boyd_2011_DOS}.
Specifically, ADMM splits the original optimization problem into reduced subproblems that are solved in sequence and whose solution is efficient and simple to implement. Furthermore, ADMM naturally handles non-smooth terms (such as the $1$-norm term in \eqref{main_con}) by employing proximal methods. {Very recent research has focussed on applying ADMM to nonlinear, nonconvex equality-constrained optimization problems, similar to problem formulation (\ref{main_con}). For instance, \cite{Wang2021} analyzes the convergence of a nonlinear objective function subject to a 2-norm squared equality constraint; because of our approach in handling each subproblem (see Section \ref{TTV}), such analysis can be adapted to our problem formulation \eqref{main_con}.}

The value of the parameter $\beta$ has a great impact on the 
quality of the regularized solution for the problem at hand. The authors of \cite{Gholami_2013_BCT} 
proposed to determine $\beta$ via L-curve analysis, still assuming that an estimate of the error bound $\varepsilon$ is known a priori. Specifically, an unconstrained counterpart of \eqref{main_con} is considered, which involves two parameters (one for weighting the fit-to-data and the regularization terms, and one for balancing the two regularization terms); the problem is solved by applying both the discrepancy principle and the L-curve (taking a range of values assigned to $\beta$), requiring the solution of many instances of the considered optimization problem (one for each value of $\beta$) to generate the L-curve. 
Recently, the authors of 
\cite{huska2019convex}
presented a procedure for tuning $\beta$ in image decomposition problems, i.e., problems like \eqref{main_con} where $\bold{G}$ is the identity matrix. However, their method assumes that a good estimate of the minimum non-zero gradient norm of the blocky component $\bold{m}_1$ is known a priori. 
Therefore, their method is limited to problems for which this quantity can be estimated directly from the input data.

In this paper, we present an extremely simple and efficient strategy for automatic selection of $\beta$ using robust statistics. Our basic idea is that, since the desired model $\bold{m}$ consists of the blocky component $\bold{m}_1$ and the smooth component $\bold{m}_2$, from a statistical point of view the model gradient 
$\bold{D}_1 \bold{m}$ 
is from a mixture of non-Gaussian and Gaussian distributions. Specifically, the entries of 
$\bold{D}_1 \bold{m}_1$ 
are considered as anomalies/outliers in the model gradient and therefore statistical anomaly detection tools are used to identify them. The robust z-score \cite{iglewicz1993detect} allows us to optimally determine the lower limit of the anomalous entries (or the upper limit of normal/Gaussian distributed entries) in the gradient vector. In particular, the optimum $\beta$ is determined such that 
$\|\bold{D}_1\bold{m}_2\|_{\infty}$
is equal to the minimum value of the anomalies (or, equivalently, the maximum value of the normal entries) determined by the z-score. 
A simple algorithm is presented to achieve this, which updates the value of $\beta$ at each ADMM iteration, leading to an optimal separation of the model components. Extensive numerical examples from imaging, tomography and geophysical inverse problems are presented to show the excellent performance of the proposed method for solving ill-posed inverse problems with piecewise smooth solutions.

The rest of the paper is organized as follows. In Section \ref{TTV}, we present more details about the formulation and ADMM implementation of the Tikhonov-TV regularization method \eqref{main_con}. In Section \ref{ZSCORE}, we introduce the new technique for the selection of the balancing parameter $\beta$. 
Experimental results are displayed in Section \ref{NUMERICAL}. Finally, some concluding remarks are proposed in Section \ref{CONS}.

\section{Solving the Tikhonov-TV regularized problem using ADMM} \label{TTV}

{In this section we provide some details about the formulation and the solution of problem \eqref{main_con}, where the parameter $\beta$ is fixed.} We start by precisely defining the matrices $\bold{D}_{1}$ and $\bold{D}_{2}$ introduced in Section \ref{intro}. For problems (\ref{main_eq}) formulated in 1D, we take

\begin{equation}  \label{D}
\bold{D}_{1}:=\bold{D}_{1,\oneD}^{(N)}:=
\begin{pmatrix}
-1 & 1 & & & \\
& -1 & 1 & & \\
& & \ddots & \ddots & \\
& & & -1 & 1 \\
& & &  & 0
\end{pmatrix}\in \mathbb{R}^{N\times N}
\end{equation}
and
\begin{equation} \label{Delta}
\bold{D}_{2}:=\bold{D}_{2,\oneD}^{(N)}:=
\underbrace{(\bold{D}_{1,\oneD}^{(N)})}_{=:\bar{\bold{D}}_{1}}\underbrace{\bold{D}_{1,\oneD}^{(N)}}_{=\bold{D}_1}= 
\left(
\begin{array}{rrrrrr}
1 & -2 & 1 & & &\\
& 1 & -2 & 1& &\\
& & \ddots & \ddots & \ddots &  \\
& &  & 1 & -2 & 1\\
& & &  & 1 & -1\\
& & &  &  &0
\end{array}\right)
\in \mathbb{R}^{N\times N},
\end{equation}
where the alternative notations $\bold{D}_{1,\oneD}^{(N)}$ and $\bold{D}_{2,\oneD}^{(N)}$ highlight the differentiation order, the dimensionality of the problem, and the size of the matrices. Although, for the 1D case, $\bar{\bold{D}}_{1}=\bold{D}_1=\bold{D}_{1,\oneD}^{(N)}$, we keep these notations to better match the 2D case (see \eqref{Delta2}). Note that $\bold{D}_{1}$ and $\bold{D}_{2}$ as expressed above are associated to a rescaled forward first order and second order finite difference scheme, respectively. For problems (\ref{main_eq}) formulated in 2D, denoting Kronecker products by $\otimes$ and identity matrices of size $N$ by $\bold{I}_{N}$, we take
\begin{equation}
\bold{D}_{1}:=
\begin{pmatrix}
\bold{D}_{1,\oneD}^{(N_x)}\;&\!\!\!\otimes\!\!\!\!& \bold{I}_{N_z}\\
\bold{I}_{N_x}&\!\!\otimes\!\!\!&\,\bold{D}_{1,\oneD}^{(N_z)}
\end{pmatrix}\in \mathbb{R}^{2N\times N},
\end{equation}
which represents the discrete partial first order derivatives in the $x$ (or horizontal) direction (matrix on the top) and in the $z$ (or vertical) direction (matrix on the bottom). Similarly to the 1D case, we take
\begin{align}
\bold{D}_{2}&:=
\begin{pmatrix}
(\bold{D}_{1,\oneD}^{(N_x)})\bold{D}_{1,\oneD}^{(N_x)}&\!\!\otimes\!\!& \bold{I}_{N_z}\\
\bold{I}_{N_x}&\!\!\otimes\!\!&(\bold{D}_{1,\oneD}^{(N_z)})\bold{D}_{1,\oneD}^{(N_z)}
\end{pmatrix}\nonumber\\
&\;=
\underbrace{\begin{pmatrix}
(\bold{D}_{1,\oneD}^{(N_x)})\otimes \bold{I}_{N_z} & \bold{0}\\
\bold{0} & \bold{I}_{N_x}\otimes (\bold{D}_{1,\oneD}^{(N_z)})
\end{pmatrix}}_{=:\bar{\bold{D}}_{1}}
\underbrace{
\begin{pmatrix}
\bold{D}_{1,\oneD}^{(N_x)}&\!\!\otimes\!\!& \bold{I}_{N_z}\\
\bold{I}_{N_x}&\!\!\otimes\!\!&\bold{D}_{1,\oneD}^{(N_z)}
\end{pmatrix}}_{=\bold{D}_1}\in \mathbb{R}^{2N\times N}. \label{Delta2}
\end{align}
Note that the one adopted above is not the only way of defining discrete differential operators of the first and second order. Indeed, one may consider alternative operators obtained varying boundary conditions and/or varying the way first-order operators are combined \cite{Hansen_2010_DIP}.

In order to compute the regularized solution $\bold{m}$, one can solve problem \eqref{main_con} for $\bold{m}_1$ and $\bold{m}_2$ separately, or reformulate the problem to be solved directly for $\bold{m}$. 
Note that, in general, the two components $\bold{m}_1$ and $\bold{m}_2$ are not unique: indeed, for any constant signal $\bold{c}$, we have that $\bold{D}_1\bold{c}=\bold{0}$ and $\bold{D}_2\bold{c}=\bold{0}$, so that, if the pair ($\bold{m}_1$, $\bold{m}_2$) is a solution of \eqref{main_con}, then $(\bold{m}_1\pm \bold{c}, \bold{m}_2 \mp \bold{c})$ is also a solution of \eqref{main_con}. However, the sum $\bold{m}=\bold{m}_1+\bold{m}_2$ is unique; see \cite{Vasin2018ASC}. 
Because of this, in this paper, we solve directly for $\bold{m}$. 

Let us 
introduce the auxiliary variables 
\begin{equation}
\bold{e}=\bold{d-Gm},~~~\bold{g}_1=\bold{D}_1\bold{m}_1,~~~\bold{g}_2=\bold{D}_1\bold{m}_2\,,
\end{equation}
so that problem \eqref{main_con} can be written in the following equivalent form
\begin{mini} 
{\bold{m},\bold{g}_1,\bold{g}_2,\bold{e}}{\|\bold{g}_1\|_1+\frac{\beta}{2}\|\bar{\bold{D}}_1 \bold{g}_2\|_2^2}
{\label{main_con2}}{}
\addConstraint {\bold{D}_1\bold{m}}{=\bold{g}_1+\bold{g}_2}
\addConstraint {\bold{G}\bold{m}+\bold{e}}{=\bold{d}}
\addConstraint {\|\bold{e}\|_2^2}{=\varepsilon}.
\end{mini}

The augmented Lagrangian function associated with problem \eqref{main_con2} is 
\begin{align}\label{AL}
\mathcal{L}(\bold{m},\bold{g}_1,\bold{g}_2,\bold{e},\widehat{\boldsymbol{\lambda}}_1,\widehat{\boldsymbol{\lambda}}_2,\widehat{{\lambda}}_3)=\nonumber\\
\|\bold{g}_1\|_1+\frac{\beta}{2}\|\bar{\bold{D}}_1 \bold{g}_2\|_2^2
- \langle\widehat{\boldsymbol{\lambda}}_1,\bold{D}_1\bold{m}-\bold{g}_1-\bold{g}_2\rangle- \langle\widehat{\boldsymbol{\lambda}}_2,\bold{G}\bold{m}+\bold{e}-\bold{d}\rangle- \langle\widehat{\lambda}_3,\|\bold{e}\|_2^2-\varepsilon\rangle\nonumber\\
+ \frac{\mu_1}{2}\|\bold{D}_1\bold{m}-\bold{g}_1-\bold{g}_2\|_2^2 +\frac{\mu_2}{2}\|\bold{G}\bold{m}+\bold{e}-\bold{d}\|_2^2+  \frac{\mu_3}{2}(\|\bold{e}\|_2^2-\varepsilon)^2,
\end{align}
where $\langle\cdot,\cdot\rangle$ denotes the canonical inner product in $\mathbb{R}^d$ ($d=1,M,2N$), 
$\widehat{\boldsymbol{\lambda}}_1,\,\widehat{\boldsymbol{\lambda}}_2,\,\widehat{\lambda}_3$ are the Lagrange multipliers, and $\mu_1,\,\mu_2,\,\mu_3>0$ are the penalty parameters; see \cite{nocedal2006numerical}.

The $k$th ADMM iteration has the form 
\begin{subequations} 
\begin{align}
\bold{m}^{k} &= \arg\min_{\bold{m}} \:\mathcal{L}(\bold{m},\bold{g}_1^{k-1},\bold{g}_2^{k-1},\bold{e}^{k-1},\widehat{\boldsymbol{\lambda}}_1^{k-1},\widehat{\boldsymbol{\lambda}}_2^{k-1},\widehat{\lambda}_3^{k-1}), \label{admm0_1}\\
\bold{g}_1^{k} &= \arg\min_{\bold{g}_1} \:\mathcal{L}(\bold{m}^{k},\bold{g}_1,\bold{g}_2^{k-1},\bold{e}^{k-1},\widehat{\boldsymbol{\lambda}}_1^{k-1},\widehat{\boldsymbol{\lambda}}_2^{k-1},\widehat{\lambda}_3^{k-1}),\\
\bold{g}_2^{k} &= \arg\min_{\bold{g}_2} \:\mathcal{L}(\bold{m}^{k},\bold{g}_1^{k},\bold{g}_2,\bold{e}^{k-1},\widehat{\boldsymbol{\lambda}}_1^{k-1},\widehat{\boldsymbol{\lambda}}_2^{k-1},\widehat{\lambda}_3^{k-1}),\\
\bold{e}^{k} &= \arg\min_{\bold{e}} \:\mathcal{L}(\bold{m}^{k},\bold{g}_1^{k},\bold{g}_2^k,\bold{e},\widehat{\boldsymbol{\lambda}}_1^{k-1},\widehat{\boldsymbol{\lambda}}_2^{k-1},\widehat{\lambda}_3^{k-1}),\label{admm0_4}\\
\widehat{\boldsymbol{\lambda}}_1^{k}& =\widehat{\boldsymbol{\lambda}}_1^{k-1} - \mu_1 (\boldsymbol{\nabla}\bold{m}^{k}-\bold{g}_1^{k}-\bold{g}_2^{k}),\\
\widehat{\boldsymbol{\lambda}}_2^{k}& =\widehat{\boldsymbol{\lambda}}_2^{k-1} - \mu_2(\bold{G}\bold{m}^{k}+\bold{e}^k-\bold{d}). \label{admm0_5}\\
\widehat{\lambda}_3^{k}& =\widehat{\lambda}_3^{k-1} - \mu_3(\|\bold{e}^k\|_2^2-\varepsilon)\,.
\end{align}
\end{subequations} 
By combining the linear and quadratic terms in the augmented Lagrangian function \eqref{AL} as
\begin{equation}
-\langle\widehat{\boldsymbol{\lambda}}_i,\bold{x}\rangle+ \frac{\mu_i}{2}\|\bold{x}\|_2^2= \frac{\mu_i}{2}\|\bold{x}-(\frac{1}{\mu_i})\widehat{\boldsymbol{\lambda}}_i \|_2^2 - \frac{\mu_i}{2}\|\widehat{\boldsymbol{\lambda}}_i\|_2^2,
\end{equation}
and by making a change of variables $\boldsymbol{\lambda}_i=(1/\mu_i)\widehat{\boldsymbol{\lambda}}_i, i=1,2,3$,  
the ADMM iteration in \eqref{admm0_1}-\eqref{admm0_5} can be written in the simpler scaled form
\begin{subequations} 
\begin{align}
\bold{m}^{k} &= \arg\min_{\bold{m}} 
\{\frac{\mu_1}{2} \|\bold{D}_1\bold{m}-\bold{g}_1^{k-1}-\bold{g}_2^{k-1} -\boldsymbol{\lambda}_1^{k-1}\|_2^2 
+\frac{\mu_2}{2}\|\bold{Gm}+\bold{e}^{k-1}-\bold{d}-\boldsymbol{\lambda}_2^{k-1}\|_2^2\}, \label{admm1_1}\\
\bold{g}_1^{k} &= \arg\min_{\bold{g}_1} \{\frac{\mu_1}{2} \|\bold{D}_1\bold{m}^k-\bold{g}_1-\bold{g}_2^{k-1} -\boldsymbol{\lambda}_1^{k-1}\|_2^2+ \|\bold{g}_1\|_1\},\label{admm1_2}\\
\bold{g}_2^{k} &= \arg\min_{\bold{g}_2} \{\frac{\mu_1}{2} \|\bold{D}_1\bold{m}^k-\bold{g}_1^{k}-\bold{g}_2 -\boldsymbol{\lambda}_1^{k-1}\|_2^2+\frac{\beta}{2}\|\bar{\bold{D}}_1 \bold{g}_2\|_2^2\}, \label{admm1_3}\\
\bold{e}^{k} &= \arg\min_{\bold{e}} \{\frac{\mu_2}{2}\|\bold{Gm}^k+\bold{e}-\bold{d} -\boldsymbol{\lambda}_2^{k-1}\|_2^2+\frac{\mu_3}{2}\left(\|\bold{e}\|_2^2-\varepsilon - \lambda_3^{k-1}\right)^2\}, \label{admm1_4}\\
\boldsymbol{\lambda}_1^{k}& =\boldsymbol{\lambda}_1^{k-1} + \bold{g}_1^{k}+\bold{g}_2^{k}-\bold{D}_1\bold{m}^{k},\label{admm1_5}\\
\boldsymbol{\lambda}_2^{k}& =\boldsymbol{\lambda}_2^{k-1} + \bold{d}- \bold{e}^k-\bold{G}\bold{m}^{k}, \label{admm1_6}\\
{\lambda}_3^{k}& ={\lambda}_3^{k-1} + \varepsilon-\|\bold{e}^k\|_2^2\,.  \label{admm1_7}
\end{align}
\end{subequations} 
In the following, we explain how the subproblems \eqref{admm1_1}-\eqref{admm1_4} can be solved efficiently. 

Regarding $\bold{m}$, since the associated subproblem \eqref{admm1_1} is a differential Tikhonov regularization problem in standard form, the optimality conditions for \eqref{admm1_1} lead to
\begin{equation} \label{m_sub}
\bold{m}^{k}\!=\!\left(\mu_1 \bold{D}_1^T\bold{D}_1 + \mu_2 \bold{G}^T\bold{G}  \right)^{-1}\!\!\!
\left(\mu_1 \bold{D}_1^T( \bold{g}_1^{k-1}\!\! + \bold{g}_2^{k-1}+\boldsymbol{\lambda}_1^{k-1})\! + \mu_2\bold{G}^T(\bold{d}-\bold{e}^{k-1}\!\!+\boldsymbol{\lambda}_2^{k-1})\right)\!.
\end{equation}
Although $\bold{m}^{k}$ has the closed form expression given above, in practice and in a large-scale setting, depending on the properties of $\bold{G}$, it may be necessary to apply an iterative solver (such as CG or CGLS) to compute $\bold{m}^{k}$; additional details about this will be provided in Section \ref{NUMERICAL}.  

Regarding $\bold{g}_1$, since the associated subploblem \eqref{admm1_2} is, by definition, the proximal operator associated with the $\ell_1$- norm penalty, this can be computed as follows
\begin{align} \label{p1_sub}
\bold{g}_1^{k}= \mathcal{T}_{\frac{1}{\mu_1}}(\bold{D}_1\bold{m}^{k}-\bold{g}_2^{k-1}-\boldsymbol{\lambda}_1^{k-1}),
\end{align}
where the operator $\mathcal{T}_{\frac{1}{\mu_1}}$ is the soft-thresholding operator defined component-wise as:
\begin{equation}
[\mathcal{T}_{\frac{1}{\mu_1}}(\bold{x})]_i= x_i~\max(1 - \frac{1}{\mu_1|x_i|},0).
\end{equation}

Regarding $\bold{g}_2$, again the associated subproblem \eqref{admm1_3} is differentiable, and the optimality condition for \eqref{admm1_3} leads to
\begin{equation} \label{p2_sub}
\bold{g}_2^{k}= \left(\bold{I} + (\beta/\mu_1) \bar{\bold{D}}_1^T\bar{\bold{D}}_1\right)^{-1}\left({\bold{D}_1}\bold{m}^{k}-\bold{g}_1^{k}-\boldsymbol{\lambda}_1^{k-1}\right).
\end{equation}
This is a first-order Tikhonov filter and can be calculated efficiently by different methods including those based on the fast Fourier transform or the discrete cosine transform; see \cite{SIRev99} and the references therein. 

Regarding $\bold{e}$, the associated subproblem \eqref{admm1_4} is differentiable but {\color{blue}nonconvex}, and the optimality condition of \eqref{admm1_4} guarantees that $\bold{e}^{k}$ satisfies
\begin{equation} \label{eqe}
\left(1 + \frac{2\mu_3}{\mu_2} (\|\bold{e}^{k}\|_2^2-\varepsilon - \lambda_3^{k-1})\right)\bold{e}^{k} = \bold{d}-\bold{Gm}^k+\boldsymbol{\lambda}_2^{k-1}\,.
\end{equation}
The solution to this system of equations may not be unique. However, we claim that 
\begin{equation} \label{egamma}
\bold{e}^{k}= \gamma_k(\bold{e}^{k}) (\bold{d}-\bold{Gm}^k+\boldsymbol{\lambda}_2^{k-1})\,,
\end{equation}
where the scale parameter $\gamma_k=\gamma_k(\bold{e}^{k})$ is the maximum real root of the depressed cubic equation (see \cite{poly})
\begin{equation} \label{cubic} 
\gamma^3\,+\, \frac{\mu_2-2\mu_3(\varepsilon+\lambda_3^{k-1})}{2\mu_3 \|\bold{d}-\bold{Gm}^k+\boldsymbol{\lambda}_2^{k-1}\|_2^2}\,\gamma\, + \,\frac{-\mu_2}{2\mu_3 \|\bold{d}-\bold{Gm}^k+\boldsymbol{\lambda}_2^{k-1}\|_2^2}=0\,,
\end{equation}
satisfies \eqref{eqe} and is also the global minimizer of  \eqref{admm1_4}. Indeed, expression \eqref{egamma} comes from the fact that $\bold{e}^{k}$ on the left-hand side of \eqref{eqe} is premultiplied by a real scalar dependent on $\|\bold{e}^{k}\|_2^2$. By substituting $\bold{e}^{k}$ expressed as in \eqref{egamma} into \eqref{eqe} leads to the determination of $\gamma_k$ as a root of the depressed cubic equation \eqref{cubic}. Furthermore, plugging $\bold{e}^{k}$ from \eqref{egamma} into the objective function \eqref{admm1_4} (here denoted by $f(\gamma_k)$), and letting 
\[
E_k=\|\bold{d}-\bold{Gm}^k+\boldsymbol{\lambda}_2^{k-1}\|_2^2,\quad
p_k=\frac{\mu_2-2\mu_3(\varepsilon+\lambda_3^{k-1})}
{E_k},\quad
q_k=-\frac{\mu_2}{E_k}
\]
to simplify the notations, we get that 
\begin{align}
f(\gamma_k)&=\frac{\mu_2}{2}\|\bold{Gm}^k+\gamma_k (\bold{d}-\bold{Gm}^k +\boldsymbol{\lambda}_2^{k-1}) -\bold{d} -\boldsymbol{\lambda}_2^{k-1}\|_2^2\nonumber\\
&\quad+\frac{\mu_3}{2}\left(\|\gamma_k (\bold{d}-\bold{Gm}^k +\boldsymbol{\lambda}_2^{k-1})\|_2^2-\varepsilon - \lambda_3^{k-1}\right)^2\nonumber\\ 
&= \frac{\mu_2 E_k}{2}(\gamma_k-1)^2 + \frac{\mu_3}{2}(E_k\gamma_k^2-\varepsilon - \lambda_3^{k-1})^2\nonumber\\
&=\mu_3 E_k^2 \gamma_k \left(\frac{1}{2}\gamma_k^3 +  \frac{\mu_2-2\mu_3(\varepsilon+\lambda_3^{k-1})}{2\mu_3 E_k}\gamma_k - \frac{\mu_2}{\mu_3 E_k} \right)
+ \frac{1}{2}(\mu_2 E_k + \mu_3(\varepsilon + \lambda_3^{k-1})^2) \nonumber\\
&=\mu_3 E_k^2 \gamma_k \left(-\frac{1}{2}\gamma_k^3 + \gamma_k^3 +  p_k\gamma_k + q_k - \frac{\mu_2}{2\mu_3 E_k} \right)
+ \frac{1}{2}(\mu_2 E_k + \mu_3(\varepsilon + \lambda_3^{k-1})^2) \nonumber\\
&=-\mu_3 E_k^2 \gamma_k \left(\frac{1}{2}\gamma_k^3 + \frac{\mu_2}{2\mu_3 E_k} \right)
+ \frac{1}{2}(\mu_2 E_k + \mu_3(\varepsilon + \lambda_3^{k-1})^2)\,, 
\end{align}
where we have used the fact that $\gamma_k$ solves \eqref{cubic}. It follows that $f(\gamma_k)$ is minimized when evaluated at the largest root of \eqref{cubic}. 

\section{Balancing Parameter Selection by Robust Statistics} \label{ZSCORE}
The Tikhonov-TV regularized problem \eqref{main_con} introduced in Section \ref{intro} stems from the statistical assumption that the model gradient $\bold{g}=\bold{D}_1\bold{m}$ is a mixture of two components resulting from different distributions: a non-Gaussian distributed (sparse) component $\bold{g}_1=\bold{D}_1\bold{m}_1$ and a  Gaussian distributed (non-sparse) component $\bold{g}_2=\bold{D}_1\bold{m}_2$. We determine the balancing parameter $\beta$ by using statistical tools that allow to optimally separate these two components of the gradient. Specifically, we assume that the non-zero entries of $\bold{g}_1$, which are associated with jumps in the regularized solution, can be considered as anomalies (or outliers) in the gradient vector $\bold{g}$. The lower {value} of these anomalous entries can be optimally determined by robust statistics \cite{rousseeuw2018anomaly}. 

According to the classical z-score statistic, an element of a data set is considered an anomaly if it falls outside a distance (e.g., 2.5 standard deviations) from the mean. A smaller distance (e.g., 2.0 standard deviations) can be used if the model size is small and a larger distance (e.g., 3.0 or 3.5 standard deviations) can be used if the model size is large. Therefore, its value can be determined  according to the user's perspective; {through this paper we denote this threshold distance by ${\tau_{\tiny{\text{nrm}}}}$ and, unless otherwise stated, we set its value to 2.5} for all numerical examples.
A main difficulty in anomaly detection is that the z-score is sensitive to extreme {values} because the mean and standard deviation are sensitive to extreme values.

The robust z-score \cite{iglewicz1993detect}  is based on robust estimators rather than the mean and the standard deviation. 
A robust measure of mean is the median. For a given vector $\bold{g}$, the median of $\bold{g}$, denoted by $\text{median}(\bold{g})$, is the value such that  at least half of the entries in $\bold{g}$ are smaller than or equal to $\text{median}(\bold{g})$, and that at least half of the entries in $\bold{g}$ are larger than or equal to $\text{median}(\bold{g})$. 
Thus, if $\bold{g}$ has an odd number of sorted elements, then $\text{median}(\bold{g})$ is the middle value; otherwise it can be the average of the two middle values. 
A robust measure of standard deviation is known as the median absolute deviation (MAD). For a given vector $\bold{g}$, the MAD of $\bold{g}$, denoted by $\text{MAD}(\bold{g})$, is given by
\begin{equation}
\text{MAD}(\bold{g})=1.4826~ \text{median}(\bold{g} - \text{median}(\bold{g}))\,,
\end{equation}
where the constant value 1.4826 is used to make the estimator consistent for Gaussian distributions; see \cite{rousseeuw1993alternatives}. 
Other alternatives to the MAD can also be used for estimation of the standard deviation; see, again, \cite{rousseeuw1993alternatives}. 
%
%
Using these robust measures of the mean and standard deviation, a robust z-score is assigned to each gradient sample, calculated as 
\begin{equation}\label{rob_z-score}
z_i=\frac{[\bold{g}]_i - \text{median}(\bold{g})}{\text{MAD}(\bold{g})}\,,
\end{equation}
and every element for which $|z_i|$ is larger than the threshold value ${\tau_{\tiny{\text{nrm}}}}$ is considered as anomaly. Consequently, we define $\text{nrm}(\bold{g})$ as the vector comprising the `normal' (i.e., the `not-anomalous') entries of $\bold{g}$, i.e.,
\begin{equation} \label{normal_set}
\text{nrm}(\bold{g})=\{[\bold{g}]_i:~ |z_i| \leq {\tau_{\tiny{\text{nrm}}}}\}.
\end{equation}

According to the insight given at the beginning of this section, the optimal balancing parameter 
to be used in \eqref{main_con} should be such that the smooth components $\bold{g}_2$ of the computed solution gradient $\bold{g}$ are classified as `normal' \eqref{normal_set} according to the robust z-score \eqref{rob_z-score}. We impose this condition by defining $\beta_{\opt}$ as a root of the equation
\begin{equation} \label{root}
0=\|\bold{g}_2(\beta)\|_{\infty} - \|\text{nrm}(\bold{g}(\beta))\|_{\infty}=:\phi(\beta)\,,\quad\mbox{i.e.,}\quad \phi(\beta_{\opt})=0\,.
\end{equation}
In the above equation, the notations $\bold{g}_2(\beta)$ and $\bold{g}(\beta)$ have been employed to highlight the dependence of the vectors $\bold{g}_2$ and $\bold{g}$ on $\beta$; in the following, both these notations will be freely used depending on the context. In \eqref{root}, the vectors $\bg_2$ and $\nrmg$ are evaluated in the max norm to enforce that the value of the maximum entry in the smooth and normal components of $\bg$ coincide. 

{
\begin{proposition}
The function $\phi(\beta)$ defined in \eqref{root} has at least one positive root. 
\end{proposition}} 
\begin{proof}
 Once the ADMM iterations \eqref{admm1_1}-\eqref{admm1_7} have converged, {we} can write the computed $\bold{g}_1$ and $\bold{g}_2$ as 
\begin{eqnarray}
\bold{g}_1&=& \mathcal{T}_{\frac{1}{\mu_1}}(\bold{D}_1\bold{m}-\bold{g}_2-\boldsymbol{\lambda}_1)
\label{eq:g1_1}\\
\bold{g}_2&=& \left(\bold{I} + (\beta/\mu_1) \bar{\bold{D}}_1^T\bar{\bold{D}}_1\right)^{-1}\left({\bold{D}_1}\bold{m}-\bold{g}_1-\boldsymbol{\lambda}_1\right)
\label{eq:g2_1}
\end{eqnarray}
where  the updates rules \eqref{p1_sub} and \eqref{p2_sub} 
have been exploited. 
From the above equalities it can be clearly seen  that, when $\beta= 0$, $\bold{g}_2(0)=\bold{D}_1\bm - \bg_1-\blam_1$, so that 
\[
\bg_1 = \mathcal{T}_{\frac{1}{\mu_1}}(\bold{D}_1\bold{m}-\bold{g}_2-\mathbf{\lambda}_1)=\mathcal{T}_{\frac{1}{\mu_1}}(\bold{g}_1)\,,
\] 
which holds if and only if $\bg_1=\mathbf{0}$. As a consequence, $\bold{g}(0)=\bold{g}_2(0)$. After defining the diagonal matrix $\bold{D}(\bold{g}_2)$ whose diagonal entries are 
\[
[\bold{D}(\bold{g}_2)]_{i,i}=\begin{cases}1 & \mbox{if $[\bold{g}_2]_i\in\text{norm}(\bold{g})$}\\
0 & \mbox{otherwise}
\end{cases}\,,
\]
and using standard norm inequalities, we get
\[
\|\text{nrm}(\bold{g}(0))\|_{\infty}=\|\bold{D}(\bold{g}_2)\bold{g}_2(0)\|_{\infty}\leq \|\bold{D}(\bold{g}_2)\|_{\infty}\|\bold{g}_2(0)\|_{\infty} =\|\bold{g}_2(0)\|_{\infty}\,,
\]
so that
\[
\phi(0)=\|\bold{g}_2(0)\|_{\infty}-\|\text{nrm}(\bold{g}_2(0))\|_{\infty}\geq 0\,.
\] 
Let us denote by $\bv_i$ and $\sigma_i$, $i=1,\dots,N$ the eigenvectors and nonnegative eigenvalues of $\bar{\bold{D}}_1^T\bar{\bold{D}}_1$, respectively. It follows that $\bg_2$ can be expressed as 
\[
\bg_2=\sum_{i=1}^N\left(1+\frac{\beta}{\mu_1}\sigma_i\right)^{-1}(\bv_i^T(\bD_1\bm-\bg_1-\blam_1))\bv_i,
\]
so that, when $\beta\rightarrow +\infty$, $\bold{g}_2(\beta)\rightarrow \mathbf{0}$ and $\|\bold{g}_2(\beta)\|_{\infty}\rightarrow 0$. As a consequence, 
\[
\phi(\beta)=\|\bold{g}_2(\beta)\|_{\infty}-\|\text{nrm}(\bold{g}(\beta))\|_{\infty}\rightarrow-\|\text{nrm}(\bold{g}(\beta))\|_{\infty}\leq 0\quad\mbox{as}\quad\beta\rightarrow +\infty\,.
\] 
The result follows from the continuity of $\phi(\beta)$, applying the intermediate value theorem. 
\end{proof}
Under the reasonable assumption that $\|\nrm(\bg(\beta))\|_{\infty}$ is an increasing function of $\beta$, it is possible to prove that $\phi(\beta)$ defined in \eqref{root} is a strictly decreasing function of $\beta>0$ and, therefore, the zero of $\phi(\beta)$ is unique; see the arguments in Appendix \ref{append} for more details. 

In order to solve \eqref{root}, we first reformulate it as an equivalent fixed point problem of the form 
\begin{equation}\label{fixedPform}
\beta\left(\|\text{nrm}(\bold{g}(\beta))\|_{\infty}  +\|\bold{g}_2(\beta)\|_{\infty}\right) =\beta\left(2\|\bold{g}_2(\beta)\|_{\infty}\right)\,.
\end{equation}
Then, starting from an initial value $\beta^0>0$, we solve \eqref{fixedPform} via an iteration of the form
\begin{equation}\label{fixedPit}
\beta^{j+1}=\frac{1}{2}\beta^j + \frac{1}{2}\underbrace{\left(\frac{4\|\bg_2(\beta^j)\|_{\infty}}{\|\bg_2(\beta^j)\|_{\infty}+\|\nrm(\bg(\beta^j))\|_{\infty}}-1\right)\beta^j}_{=:\psi(\beta^j)},\quad
j=0,1,\dots\,.
\end{equation}
Since the above update rule can be regarded as an averaged iteration algorithm and since, under mild assumptions, it can be proved that the function $\psi(\beta)$ is not expansive (see Appendix \ref{append}), thanks to the theory presented in \cite{fixedP} the iteration \eqref{fixedPit} is guaranteed to globally converge to a fixed point of \eqref{fixedPform} and, therefore, to a solution to \eqref{root}. 

Within the ADMM scheme \eqref{admm1_1}-\eqref{admm1_5} we dynamically tune the value of $\beta$ according to rule \eqref{fixedPit}, in  such a way that condition \eqref{root} is (approximately) satisfied for the final regularized solution (i.e., at the end of the ADMM iterations). Namely, suppose that, at iteration $k$, we have access to the balancing parameter $\beta^{k-1}$ and to the quantities $\bold{g}_1^{k-1}$, $\bold{g}_2^{k-1}$, $\bold{e}^{k-1}$, $\boldsymbol{\lambda}_1^{k-1}$, $\boldsymbol{\lambda}_2^{k-1}$. Then: $\bold{m}^{k}$ is computed from $\bold{g}_1^{k-1}$, $\bold{g}_2^{k-1}$, $\bold{e}^{k-1}$, $\boldsymbol{\lambda}_1^{k-1}$, $\boldsymbol{\lambda}_2^{k-1}$; $\bold{g}_1^{k}$ is computed from $\bold{m}^{k}$, $\bold{g}_2^{k-1}$, $\boldsymbol{\lambda}_1^{k-1}$; $\bold{g}_2^{k}$ is computed from $\bold{m}^{k}$, $\bold{g}_1^{k}$, $\boldsymbol{\lambda}_1^{k-1}$ and $\beta^{k-1}$.
We set the value of $\beta^{k}$ by applying one iteration of the update rule \eqref{fixedPit} and using the updated quantities available at iteration $k$ (i.e., $\bg_2^k(\beta^{k-1})$ and $\bg^k(\beta^{k-1})$). {Although a convergence proof for the ADMM scheme \eqref{admm0_1}-\eqref{admm0_5}, with $\beta$ updated according to \eqref{fixedPit} is outside the scope of this paper, numerical experiments consistently show that both ADMM applied with the fixed value $\beta=\beta_{\opt}$ obtained when \eqref{fixedPit} has converged and ADMM applied with adaptive $\beta$ selected by \eqref{fixedPit} converge to the same solution; see also Section \ref{NUMERICAL}.}  
The resulting method is summarized in Algorithm \ref{Algorithm}.

\begin{algorithm}[!h]
\vspace{.2cm}
 \caption{ADMM for solving the Tikhonov-TV regularized problem \eqref{main_con2}, with authomatic update of $\beta$.}  \label{Algorithm}
 \begin{tabbing}
Inputs: \= penalty parameters $\mu_1,\,\mu_2,\,\mu_3>0$ {for the Lagrangian}, $\varepsilon=\|\be\|_2^2$, an initial value $\beta^0$, \\
\> {threshold $\tau_{\tiny\text{nrm}}>0$ for (\ref{normal_set})}
\end{tabbing}
\vspace{-0.3cm}
 Set {$\boldsymbol{\lambda}_1^0=\bold{0}$, $\boldsymbol{\lambda}_2^0=\bold{0}$, ${\lambda}_{{3}}^0=0$}, $\bold{g}_1^0=\bold{0}$, $\bold{g}_2^0=\bold{0}$.   \\
 \For{$k= 1,2,\cdots$ \text{until a stopping criterion is satisfied} }
 {
Compute $\bold{m}^{k}$ using \eqref{m_sub}\\ 
Compute $\bold{g}_1^{k}$ using \eqref{p1_sub}\\ 
Compute $\bold{g}_2^{k}$ using \eqref{p2_sub} with $\beta=\beta^{k-1}$\\ 
{Determine $\gamma_k$ by solving \eqref{cubic} and compute $\bold{e}^{k}$ using \eqref{egamma}\\  
Update $\boldsymbol{\lambda}_1^{k} =\boldsymbol{\lambda}_1^{k-1} + \bold{g}_1^{k}+\bold{g}_2^{k}- \bold{D}_1\bold{m}^{k}$ \\
Update $\boldsymbol{\lambda}_2^{k} =\boldsymbol{\lambda}_2^{k-1} + \bold{d} - \bold{e}^k - \bold{G}\bold{m}^{k}$\\
Update ${\lambda}_3^{k} ={\lambda}_3^{k-1} + \varepsilon - \|\bold{e}^{k}\|_2^2$}\\
Update 
$\beta^{k}=\dfrac{1}{2}\beta^{k-1} + \dfrac{1}{2}\left(\dfrac{4\|\bg_2^k(\beta^{k-1})\|_{\infty}}{\|\bg_2^k(\beta^{k-1})\|_{\infty}+\|\nrm(\bg^k(\beta^{k-1}))\|_{\infty}}-1\right)\beta^{k-1}$
}
\end{algorithm}

\section{Numerical Experiments} \label{NUMERICAL}
In this section, we assess the performance of the proposed automated ADMM algorithm in solving inverse problems arising in different applications. 

Concerning the inputs of Algorithm \ref{Algorithm}, for all the experiments we set $\mu_1>\mu_2$ and $\mu_2=\mu_3= 1$. Moreover, we assume that the value $\varepsilon=\|\be\|_2^2$ is known in all but one experiments, where we test the robustness of the proposed solver with respect to inaccuracies in the value of $\varepsilon$; {similarly, we take $\tau_{\text{\tiny{nrm}}}=2.5$ in all but one experiments, where we test the robustness of the proposed solver with respect to the choice of $\tau_{\text{\tiny{nrm}}}$.} The iterations are stopped either when a maximum number of iterations is reached 
or when the computed $\bm$ stabilizes, i.e., $\nicefrac{\|\bold{m}^{k+1}-\bold{m}^k\|_2}{\|\bold{m}^k\|_2} <10^{{-4}}$. 

We compare the results obtained by applying the new method to the results obtained applying TV regularization and Tikhonov regularization (with a regularization term of the form $\|\bar{\bold{D}}_1\bold{D}_1\bm\|_2^2$). Both TV and Tikhonov regularization are implemented through ADMM  and can be regarded as special cases of Algorithm \ref{Algorithm}. Specifically, to recover the TV formulation, we set $\bg_2=\mathbf{0}$, so that minimization in \eqref{main_con2} happens on $\bm$, $\bg_1$ and $\be$ only, and the first constraint reduces to $\bold{D}_1\bm=\bg_1$; as a consequence, the updates \eqref{m_sub}, \eqref{p1_sub}, \eqref{eqe} and \eqref{admm1_5}-\eqref{admm1_7} are performed taking $\bg_2^{k-1}=\bg_2^k=\mathbf{0}$, and update \eqref{p2_sub} is discarded. Similarly, to recover the Tikhonov formulation, we set $\bg_1=\mathbf{0}$, so that minimization in \eqref{main_con2} happens on $\bm$, $\bg_2$ and $\be$ only, and the first constraint reduces to $\bold{D}_1\bm=\bg_2$; as a consequence, the updates \eqref{m_sub}, \eqref{p2_sub}, \eqref{eqe} and \eqref{admm1_5}-\eqref{admm1_7} are performed taking $\bg_1^{k-1}=\bg_1^k=\mathbf{0}$, and update \eqref{p1_sub} is discarded. To quantitatively evaluate the accuracy of each algorithm we compute the 2-norm relative error at the $k$th iteration, $k=1,2,\dots$, defined as $\nicefrac{\|\bm^k-\bm\|_2}{\|\bm\|_2}$. 

\subsection{Compressed sensing and denoising}
\begin{figure}[htb!] 
\begin{center}
\includegraphics[scale=0.8,trim=0cm 0cm 0cm 0cm, clip]{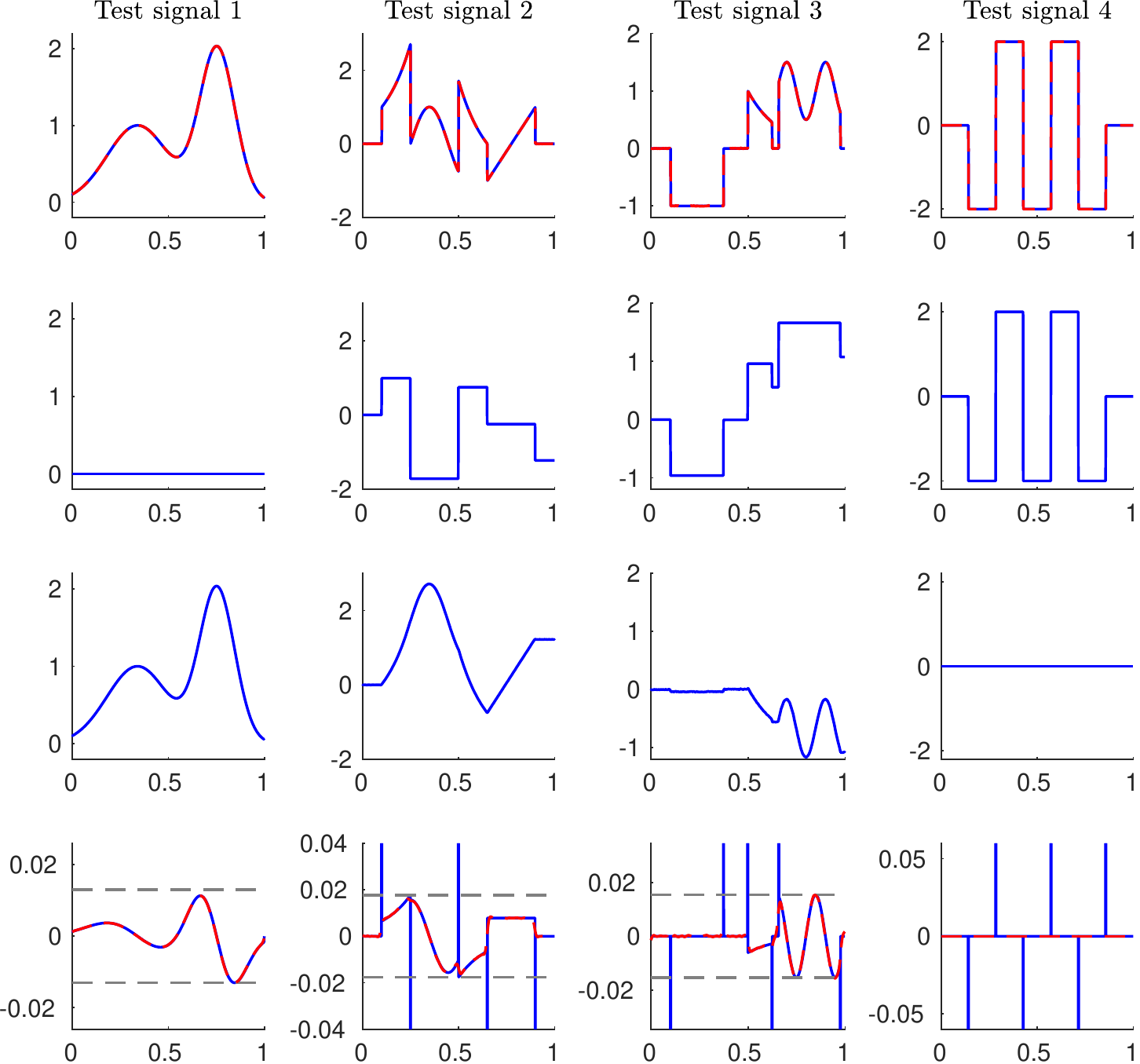}
\caption{Compressed sensing problem. First row: reconstructed signal (in blue) overlaid by the ground truth (in red). Second and third row: the blocky component $\bold{m}_1$ and the smooth component $\bold{m}_2$, respectively. Last row: the gradient of the reconstruction (in blue) overlaid by the gradient of the smooth component $\bold{m}_2$ (in red); the horizontal dashed lines show  the lower {value} determined for anomalous components.}
\label{fig:sigs_cs_signal}
\end{center}
\end{figure}

\begin{figure}[htb!] 
\begin{center}
\includegraphics[scale=0.75,trim=0cm 0cm 0cm 0cm, clip]{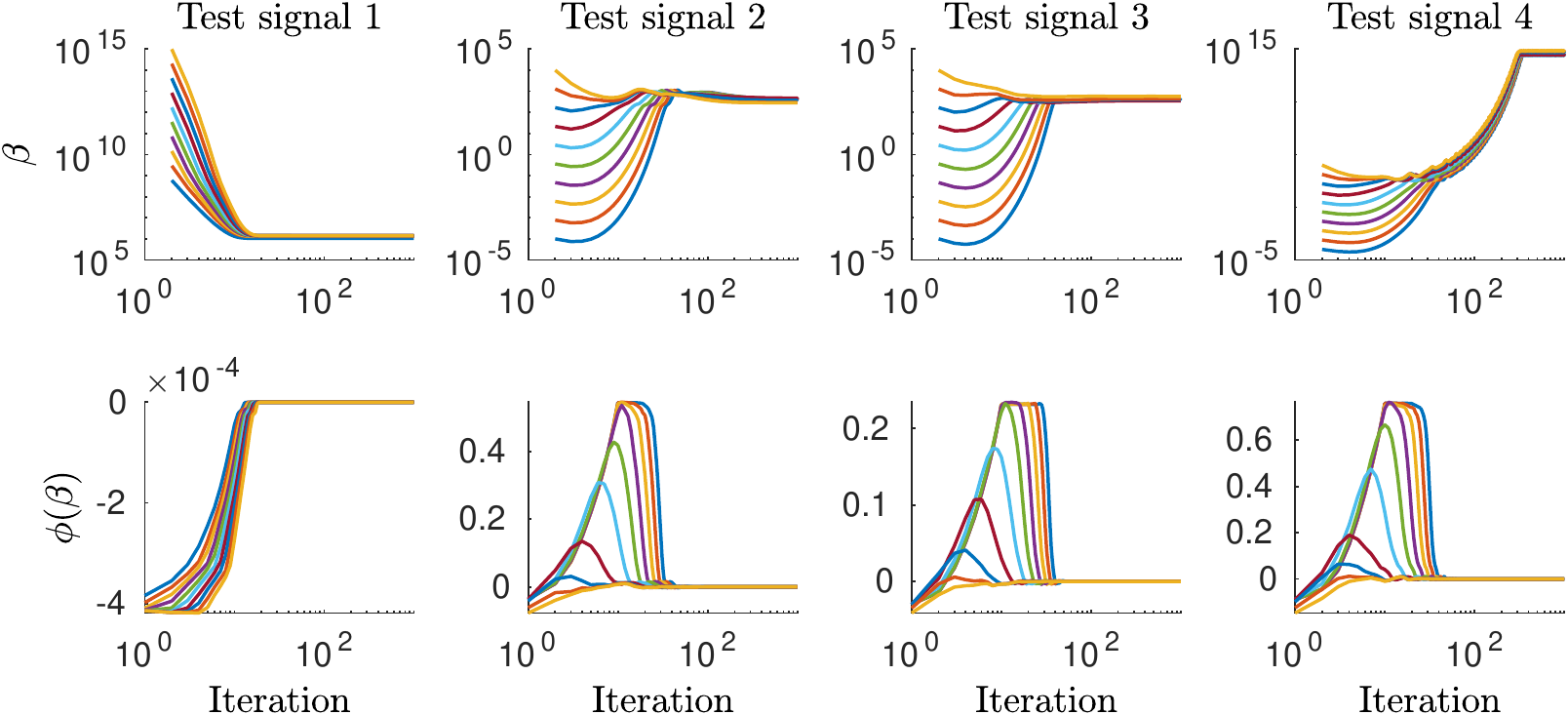}
\caption{Compressed sensing problem. Evolution of the balancing parameter $\beta$ (top row) and $\phi(\beta)$ (bottom row) versus iteration count, for different initial values $\beta_0$ of $\beta$, corresponding to the signals shown in Fig. \ref{fig:sigs_cs_signal}.
\vspace{1cm}}
\label{fig:sigs_cs_param}
\end{center}
\end{figure}

The compressed sensing theory allows to recover a sparse signal from a small number of random projections \cite{Donoho_2006_CS}. 
In this subsection, we first show the performance of the new automated Tikhonov-TV regularization method in recovering 1D signals with different features, in the framework of the compressed sensing theory.
We consider four test signals, ranging from very smooth to very rough: these are shown at the top row of Fig. \ref{fig:sigs_cs_signal} (dashed red lines). The observation matrix $\bold{G}$ consists of 1024 random column vectors of length 250, all drawn from a standard Gaussian distribution and then normalized. Each observation vector is generated as $\bold{d=Gm+e}$ where $\bold{m}$ is the test signal and $\bold{e}$ is some Gaussian white noise (such that the `noise level' $\|\be\|_{2}/\|\bold{G}\bm\|_2$ is 0.1\%).
{To evaluate the long-term behavior of} Algorithm \ref{Algorithm}, {we stop the method when it has computed} 500 iterations. The first row of Fig. \ref{fig:sigs_cs_signal} shows the reconstructed signals (in blue) overlaid by the ground truth (dashed, in red), while the nonsmooth and smooth components of the regularized solution are shown in the second and third rows, respectively. 
We can observe that nearly optimal reconstructions and decompositions are obtained for all of the signals. The model gradients of the recovered signals are shown in the last row of Fig. \ref{fig:sigs_cs_signal} (in blue), with the gradient of the recovered smooth component (dashed red curve) overlaid. In each {frame} in the bottom row of Fig. \ref{fig:sigs_cs_signal}, the lower bounds for the anomalous gradient entries as determined by the z-score \eqref{rob_z-score}, \eqref{normal_set} are shown by horizontal dashed lines.
We can observe that, in all of the cases, the z-score determined meaningful bounds, which resulted in an optimal determination of $\beta$. 
Fig.~\ref{fig:sigs_cs_param} shows the evolution of $\beta$ and $\phi(\beta)$ versus iteration count for all the four signals, with different initial guesses $\beta_0$. We can observe that formula \eqref{fixedPit} stably converges to the root of $\phi$, irrespective of the initial value $\beta_0$. 
%
%

\begin{figure}[htb!] 
\begin{center}
\includegraphics[scale=0.55,trim=0cm 0cm 0cm 0cm, clip]{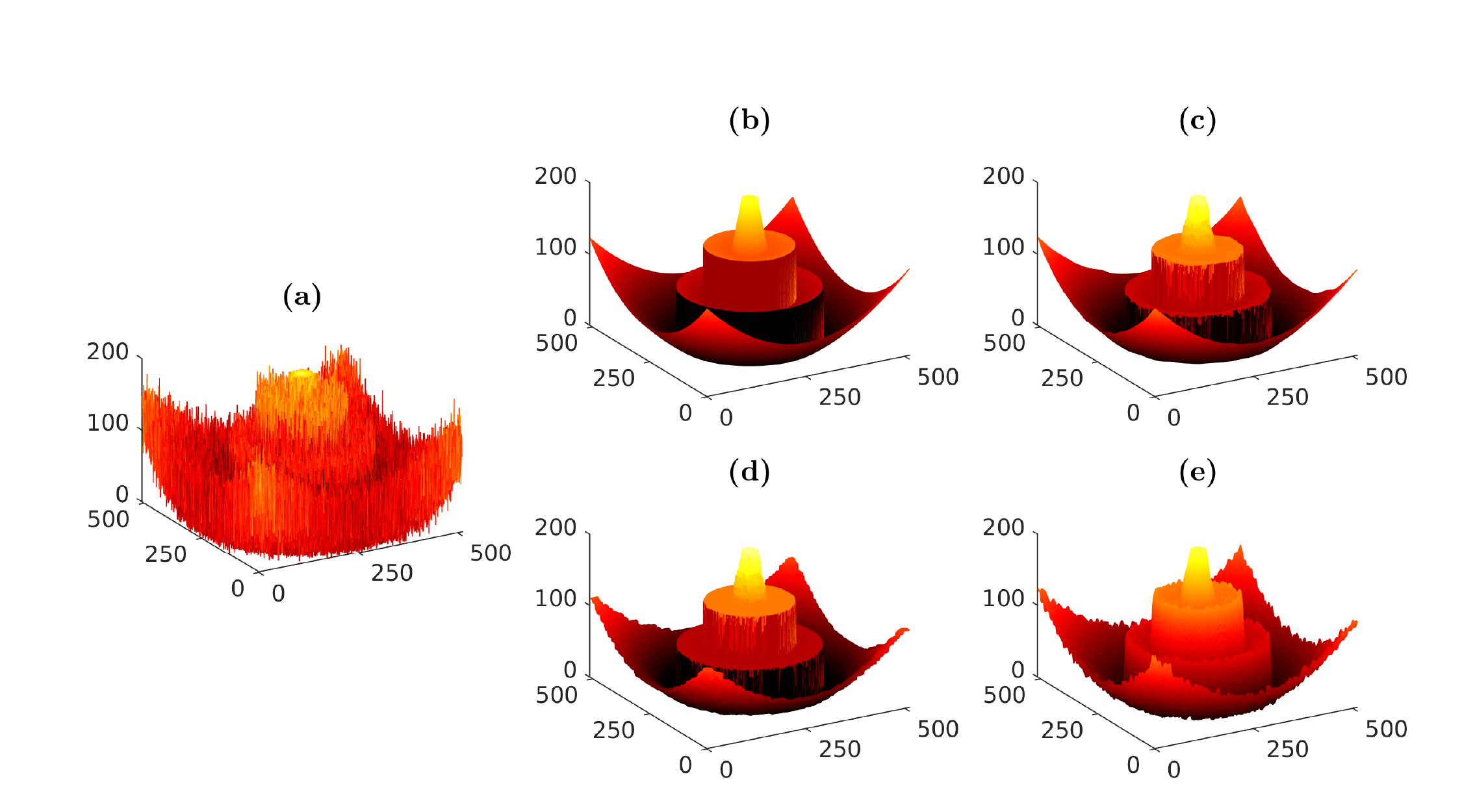}
\caption{Image denoising example. (a) degraded image by a Gaussian noise of (noise level 30\%); (b) noise-free image; (c) estimated image by the proposed method (error = $1.98\times 10^{-4}$); (d) estimated image by the TV method (error = $3.18\times 10^{-4}$); (e) estimated image by the Tikhonov method (error = $4.08\times 10^{-4}$).}
\label{fig:im_denoise}
\end{center}
\end{figure}

\begin{figure}[htb!] 
\begin{center}
\includegraphics[scale=0.58,trim=0cm 0cm 0cm 0cm, clip]{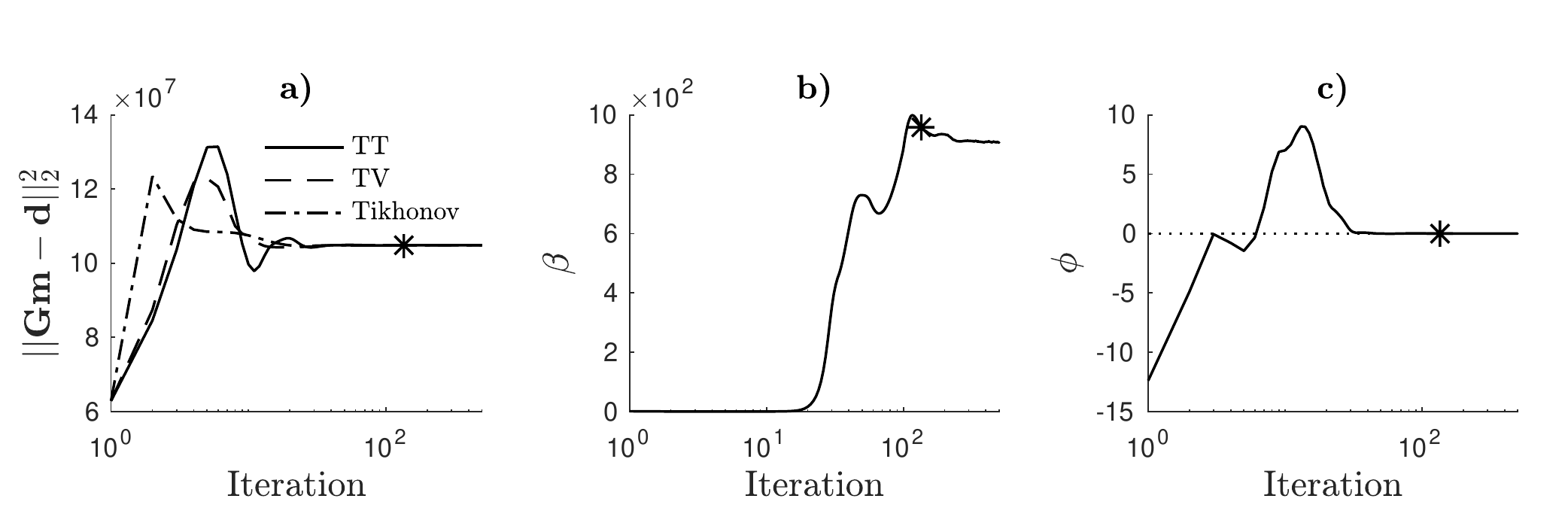}
\caption{Image denoising example. Evolution of (a) discrepancy, (b) $\beta$ and (c) $\phi(\beta)$ versus iteration for the proposed method. {The asterisk highlights the quantities computed at the 135th iteration, i.e., when the stopping criterion on the stabilization of the solution is satisfied.}}
\label{fig:im_denoise_graphs}
\end{center}
\end{figure}

\begin{figure}[htb!] 
\begin{center}
\includegraphics[scale=0.6,trim=0cm 0cm 0cm 0cm, clip]{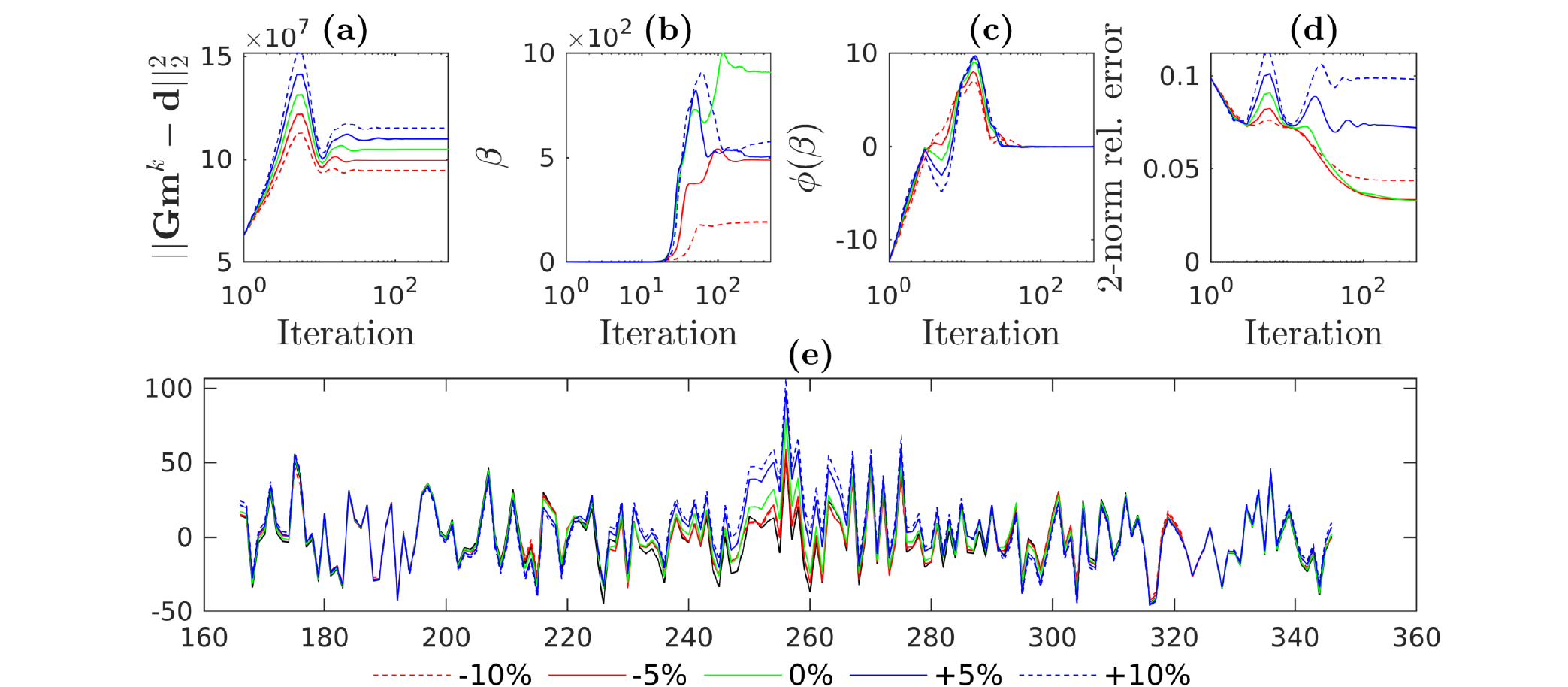}\\
\caption{Image denoising example. Sensitivity of Algorithm \ref{Algorithm} with respect to over and under estimations of $\varepsilon$. (a)-(d) evolution of (a) discrepancy, (b) $\beta$, (c) $\phi(\beta)$, and (d) 2-norm of the estimated model error. Frame (e) shows approximation of the entries 166 to 346 (at column 256) of noise vector $\be$ at the 500th (final) ADMM iteration. In {\color{blue}this} frame, the black line shows the true noise.}
\label{fig:im_denoise_sensitivity}
\end{center}
\end{figure}

\begin{figure}[htb!] 
\begin{center}
\includegraphics[scale=0.65,trim=0cm 0cm 0cm 0cm, clip]{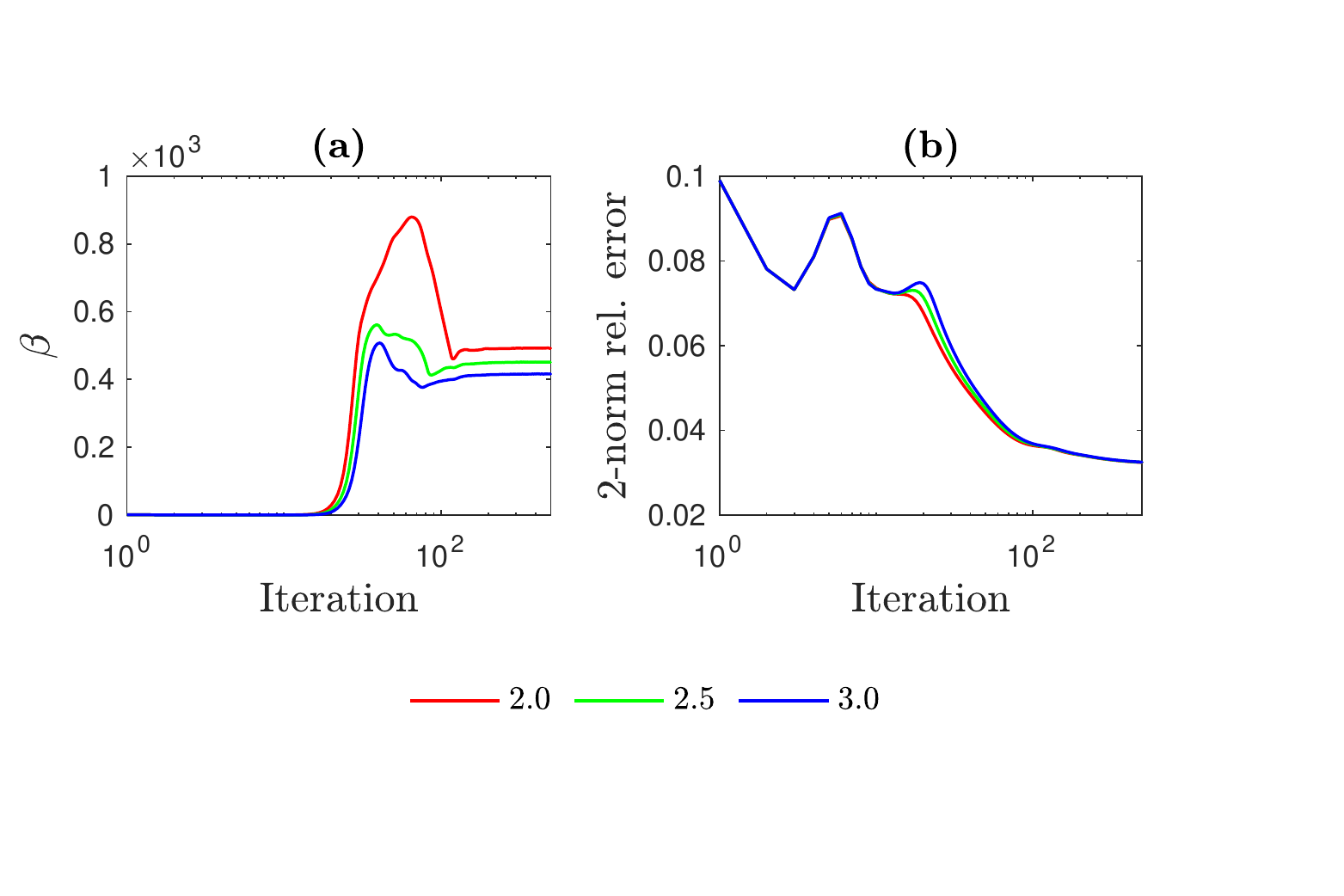}
\vspace{-1cm}
\caption{{Image denoising example. Evolution of (a) $\beta$ and (b) 2-norm relative error versus iteration for the proposed method for different values of the threshold $\tau_{\tiny{\text{nrm}}}$ in (\ref{normal_set}).}}
\label{fig:im_denoise_graphs_tau}
\end{center}
\end{figure}

\begin{figure}[htb!] 
\begin{center}
\includegraphics[scale=0.6,trim=0cm 0cm 0cm 0cm, clip]{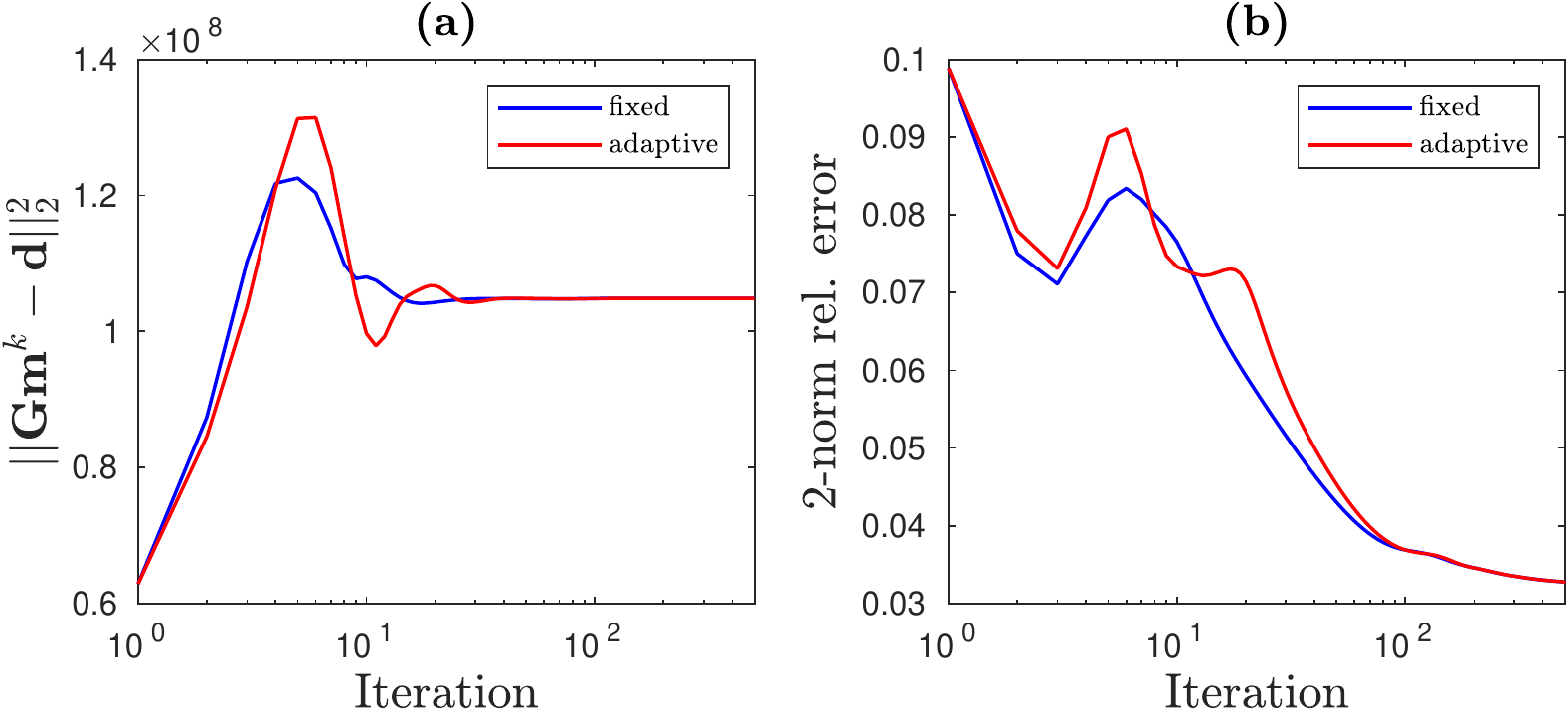}
\caption{{Image denoising example. Evolution of (a) discrepancy and (b) 2-norm of the estimated model error for Algorithm \ref{Algorithm} with iteration-dependent choice of $\beta$ and for the ADMM method \eqref{admm0_1}-\eqref{admm0_5} with fixed $\beta=\beta_{500}$.}}
\label{fig:im_denoise_beta}
\end{center}
\end{figure}

We then turn to a 2D image denoising example. The image which is to be denoised, displayed in Fig. \ref{fig:im_denoise}(a), is obtained by adding some Gaussian white noise of noise level $\|\be\|_2/\|\bm\|_2=30\%\ $ to a clean synthetic image which consists of piecewise-smooth regions, displayed in Fig. \ref{fig:im_denoise}(b). 
We use the combined Tikhonov-TV method and, to evaluate its regularization performance, we compare with the TV-only and the Tikhonov-only methods to denoise the image (as explained above). We run each method for 500 iterations {(to evaluate the long-term behavior of Algorithm \ref{Algorithm} even after the stopping criterion based on the stabilization of the solution is satisfied)} and the obtained results are depicted in Fig. \ref{fig:im_denoise}(c)-(e). As expected, the Tikhonov-TV method estimates a more accurate image by balancing between restoration of the sharp edges and of the smooth regions; instead, the TV-only and Tikhonov-only methods give more emphasis on the restoration of the sharp edges and smooth regions, respectively, and hence provide a suboptimal estimate of the noiseless image. 
The evolution of the discrepancy for all three methods is depicted in Fig. \ref{fig:im_denoise_graphs}(a). Also, Figs. \ref{fig:im_denoise_graphs}(b)-(c) show the evolution of $\beta$ and $\phi(\beta)$ versus iteration count for the Tikhonov-TV method. 
For this test problem, we experimentally study the sensitivity of {the new} Algorithm \ref{Algorithm} {with respect to a couple of algorithmic parameters that should be given in input, namely: the value of the noise magnitude $\varepsilon =\|\be\|_2^2$ and the threshold $\tau_{\text{\tiny{nrm}}}$ for the normal entries of $\bg$ (\ref{normal_set}). Concerning the former, we} consider cases where the true value of $\varepsilon$ is over- and under-estimated by 5\% and 10\%. The results of these tests are displayed in Fig.~\ref{fig:im_denoise_sensitivity}, (a)-(d). In particular, looking at frame (a), it is evident that, as the iterations proceed, the value of the squared discrepancy $\|\bold{G}\bm - \bold{d}\|_2^2$ stabilizes around the inputted values of $\varepsilon$ (leading to over- and under-fitted data when $\varepsilon$ is under- and over-estimated, respectively). Looking at frames (b)-(c) we can see that the update rule \eqref{fixedPit} for $\beta$ converges to a zero of $\phi(\beta)$ for all the considered values of $\varepsilon$, although the values of $\beta$ determined at the end of the iterations differ. The impact of an inaccurate value of $\varepsilon$ on the quality of the solution is visible in frame (d): it is clear that, for this test problem, an under-estimation of $\varepsilon$ (corresponding to an over-fitting) of the data still leads to results comparable to the case where an accurate value of $\varepsilon$ is considered; naturally, the higher the error in the estimate of $\varepsilon$, the lower the quality of the computed solution. In Fig.~\ref{fig:im_denoise_sensitivity} (e) we illustrate some of the entries of the recovered noise vector $\be$ at the end of the ADMM iterations. The new ADMM formulation allows recovery of the random noise vector $\be$ (through subproblem \eqref{admm0_4}), alongside the inverse problem solution $\bm$ and its gradient components $\bg_1$ and $\bg_2$. We can clearly see that, for this test problem, the computed approximation of $\be$ quite carefully reproduces the behavior of the unknown noise corrupting the original image. 
{Concerning the sensitivity of Algorithm \ref{Algorithm} with respect to the value of $\tau_{\text{\tiny{nrm}}}$ we consider the default value $\tau_{\text{\tiny{nrm}}}=2.5$ and $\tau_{\text{\tiny{nrm}}}=2$, $\tau_{\text{\tiny{nrm}}}=3$; the behavior of the solver for these tests are displayed in Figure~\ref{fig:im_denoise_graphs_tau}. Looking at frame (a) we can clearly see that a higher value of $\tau_{\text{\tiny{nrm}}}$ results in a lower value of the computed $\beta$ at the end of the iterations of Algorithm \ref{Algorithm}: this is expected, as a higher $\tau_{\text{\tiny{nrm}}}$ implies more entries of $\bg$ to be regarded as `normal' (see equation (\ref{normal_set})), which is achieved by imposing less penalisation on the smoothness-enforcing term in (\ref{main_con2}). Similarly, a lower $\tau_{\text{\tiny{nrm}}}$ results in a higher value of $\beta$, i.e., less entries of $\bg$ are regarded as `normal' and the smoothness-enforcing term in (\ref{main_con2}) is penalised more. Looking at frame (b) we can see that such variations in $\tau_{\text{\tiny{nrm}}}$ do not impact the quality of the solution computed by Algorithm \ref{Algorithm}. Even if not reported, a similar behavior is observed testing a wider range of $\tau_{\text{\tiny{nrm}}}$ values.} 
{Finally, we experimentally assess the convergence properties of the new ADMM scheme \eqref{admm0_1}-\eqref{admm0_5} applied with iteration-dependent choice of the balancing parameter $\beta$ according to \eqref{fixedPit}, as well as with a fixed value $\beta=\beta_{500}$, i.e., picking the value the fixed point iteration \eqref{fixedPit} converged to after 500 iterations. Fig. \ref{fig:im_denoise_beta} shows the residual and error values versus iteration number: we can clearly see that, although the two instances of ADMM display some discrepancies for approximately the first 100 iterations, they eventually converge to the same value. Even if not reported, this behavior is observed in all the considered test problems.}


\subsection{Inverse problems in Geophysics} 

We apply the new automated Tikhonov-TV regularization method for subsurface interval-velocity estimation from root-mean-square (RMS) velocities.
For a horizontally layered earth model, the RMS velocity $V(t)$ is a continuous function of time defined as 
\begin{equation} \label{Dix_c}
V(t) = \sqrt{\frac{1}{t}\int_{0}^t [v(\tau)]^2d\tau},
\end{equation}
where $v(t)$ is the instantaneous/interval velocity function and $t$ is the wave propagation time; see \cite{Dix_1955_SVS}. 
In practice, the velocity analysis of common-depth-point (CDP) gathers gives an estimate of the RMS velocity; when both the RMS velocity $V(t)$ and the interval velocity $v(t)$ are needed, the estimation of the latter from the RMS velocity is a severely ill-conditioned problem, thus proper regularization is required to stabilize the solution. Upon discretization, equation \eqref{Dix_c} reads
\begin{equation} \label{Dix_d}
V_i = \sqrt{\frac{1}{i}\sum_{j=1}^i v^2_j},
\end{equation}
where $i=1,...,N$ is the sample number and $N$ is the number of earth layers.
Squaring both sides of equation \eqref{Dix_d} results in a system of linear equations of the form $\bold{d}=\bold{G}\bold{m}$,
where $d_i=iV_i^2$, $m_j=v_j^2$ and the discrete forward operator $\bold{G}$ is a causal integration matrix, i.e., a lower triangular matrix of ones
\begin{equation} \label{G_dix}
\bold{G}=
\begin{pmatrix}
1 &   &  & & \\
1 & 1 &  & & \\
1 & 1 & 1 &  & \\
\vdots &\vdots & \ddots & \ddots & \\
1 & 1 & \cdots &  1& 1
\end{pmatrix}\in\mathbb{R}^{N\times N}.
\end{equation}

\begin{figure}[htb!] 
\begin{center}
\includegraphics[scale=0.62,trim=0cm 0cm 0cm 0cm, clip]{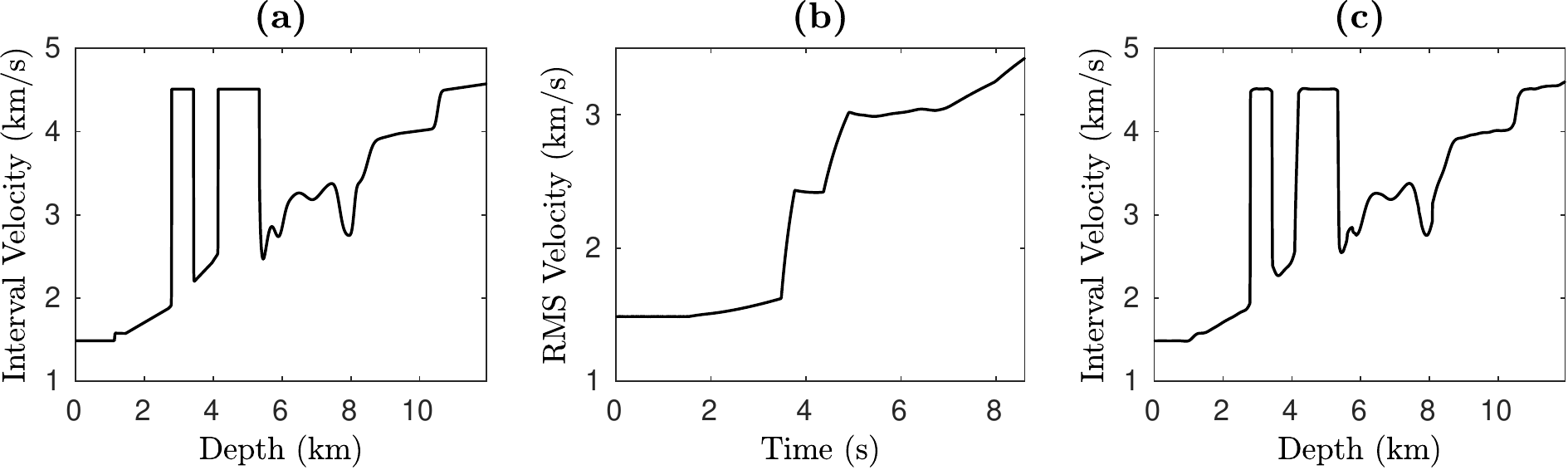}
\caption{One-dimensional Dix velocity inversion. (a) True interval velocity. (b) Decimated noisy RMS velocity. (c) Estimated interval velocity.}
\label{fig:fig_1d_CS_Dix}
\end{center}
\end{figure}

We simulate a 1D RMS velocity vector from a velocity log of the 2004 BP model \cite{Billette_2004_BPB}. The interval and RMS velocities are shown in Figs. \ref{fig:fig_1d_CS_Dix}(a) and (b), respectively. 
In order to make the simulation more realistic, only 25\% of the RMS velocity elements where used as input to the inversion, mimicking the real situations where one only picks RMS velocities at the positions of strong reflections. This leads to a compressed sensing problem for the Dix inversion \cite{Gholami_2019_3DD}. 
For this example, the forward operator is $\Phi\bold{G}$, where $\Phi$ is a 382 by 1911 sampling matrix (made up of 382 rows of an identity matrix of size 1911, picked at the locations of the velocity model with strong seismic response/reflectivity) and $\bold{G}$ is defined as in \eqref{G_dix} and has size $1911\times 1911$. 
The estimated interval velocity obtained applying 100 iterations of Algorithm \ref{Algorithm} is shown in Fig. \ref{fig:fig_1d_CS_Dix}(c); we can clearly observe that the interval velocities are estimated accurately.

\begin{figure}[htb!] 
\begin{center}
\includegraphics[scale=0.5,trim=0cm 0cm 0cm 0cm, clip]{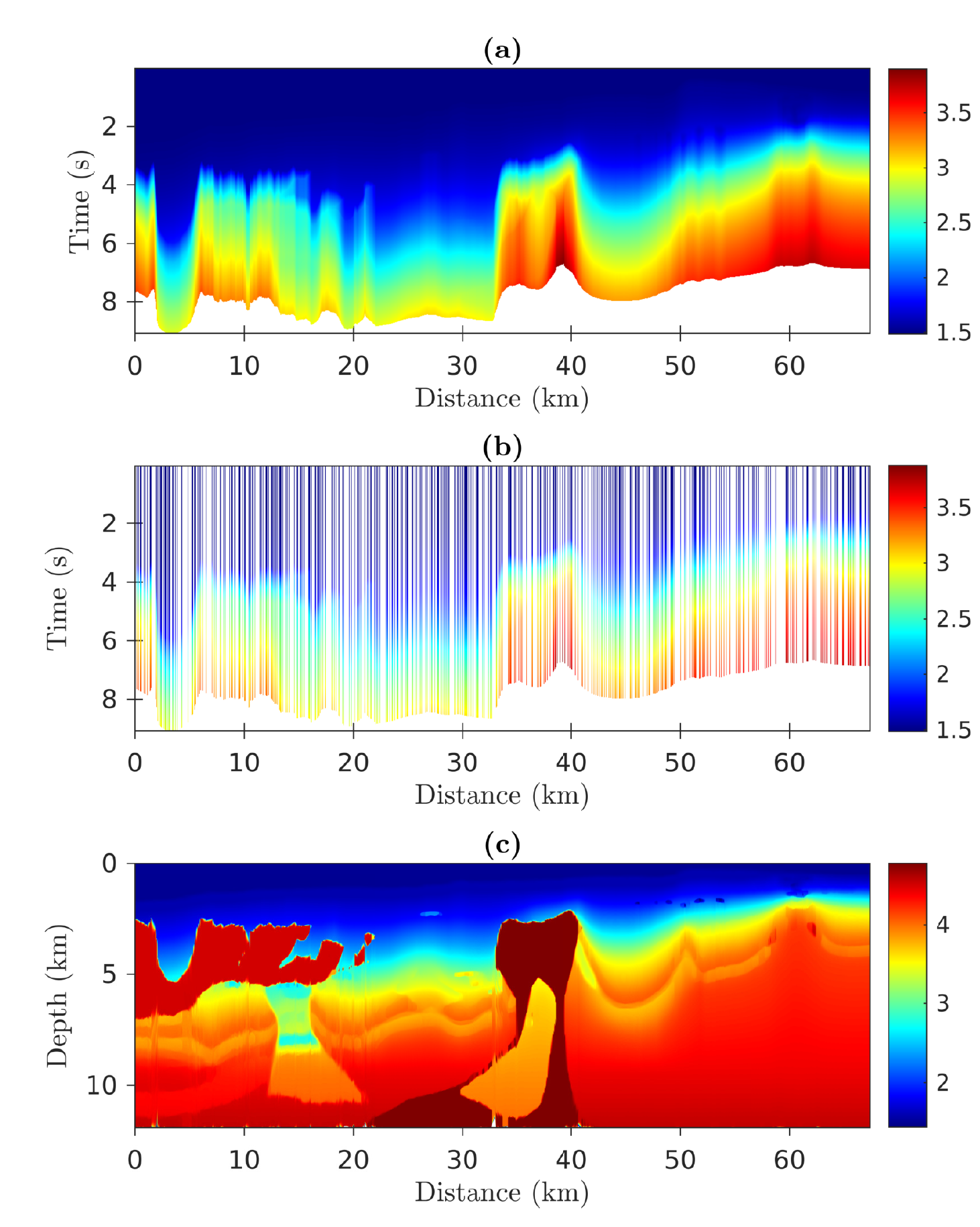}
\caption{Two-dimensional Dix velocity inversion. (a) RMS velocity field corresponding to the 2004 BP model. (b) Input RMS field from velocity scan using 30\% of CDP gathers selected at random. (c) Interval velocity field obtained by inversion using the proposed method.}
\label{fig:fig_2d_CS_Dix}
\end{center}
\end{figure}

\begin{figure}[htb!] 
\begin{center}
\includegraphics[scale=0.6,trim=0cm 0cm 0cm 0cm, clip]{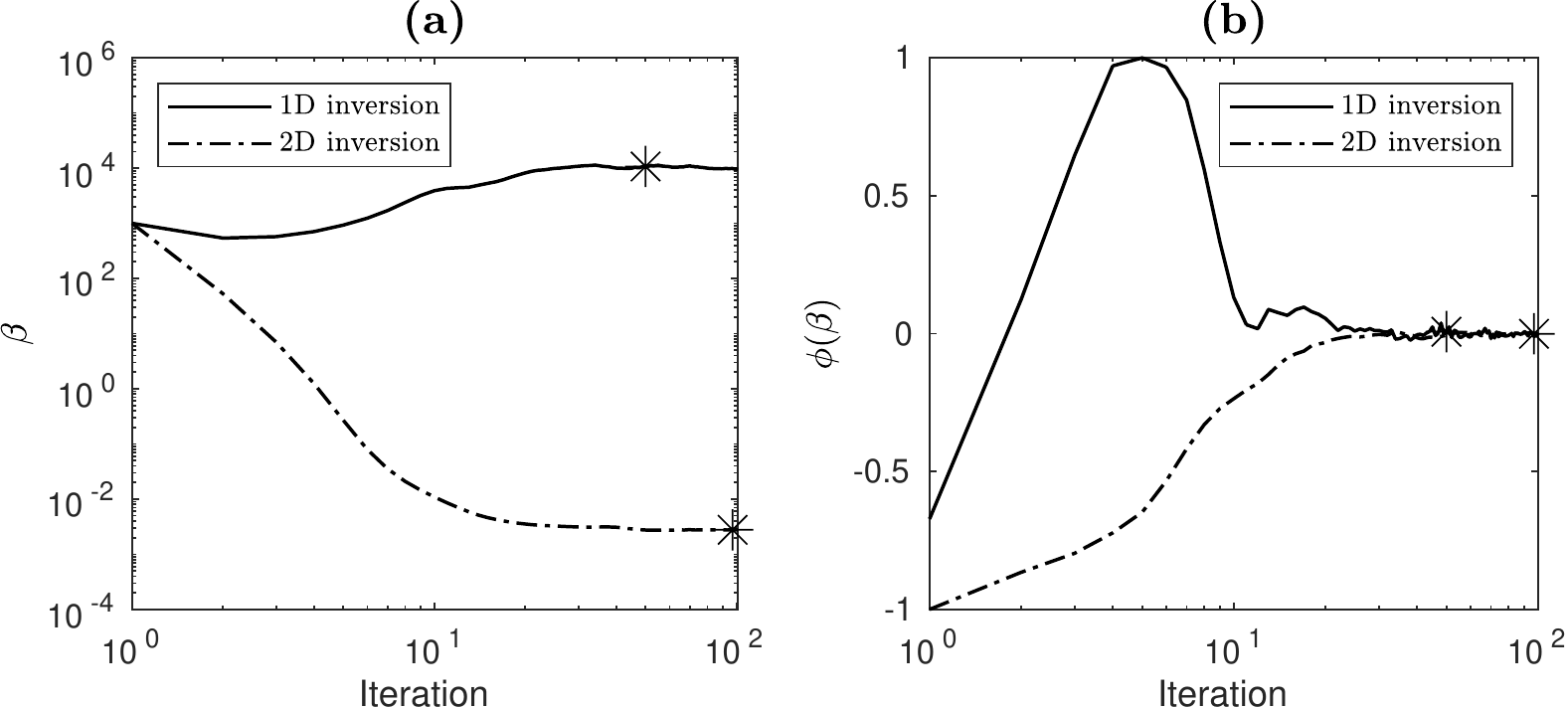}
\caption{Dix velocity inversion. Left frame: values of the balancing parameter $\beta$ versus iteration number. Right frame: values of $\phi(\beta)$ versus iteration number. Solid curves correspond to the 1D case, while dashed curves correspond to the 2D case. 
{The asterisk highlights iteration for which the stopping criterion on the stabilization of the solution is satisfied.}
}
\label{fig:fig_Dix_params}
\end{center}
\end{figure}

We then consider a 2D Dix inversion example. The full RMS velocity field is shown in Fig. \ref{fig:fig_2d_CS_Dix}(a).
We used 30\% of the traces (vertical lines, picked at random) (shown in Fig. \ref{fig:fig_2d_CS_Dix}(b)) as the input to the inversion. 
In this case, the size of the RMS velocity field is 478 $\times$ 1349 pixels and the forward operator is of the form $\Phi \otimes \bold{G}$, where $\Phi$ is a sampling matrix of size $404\times 1349$ (made up of 404 rows of an identity matrix of size 1349, picked at random) and $\bold{G}$ is defined as in \eqref{G_dix} and has size $478\times 478$.  
The estimated full interval velocity field obtained applying 100 iterations of Algorithm \ref{Algorithm} is shown in Fig. \ref{fig:fig_2d_CS_Dix}(c)
; also for this test problem we can clearly see that the interval velocity model is reconstructed accurately. 
Fig.~\ref{fig:fig_Dix_params} shows the evolution of $\beta$ and $\phi(\beta)$ versus iteration for the considered 1D and 2D Dix inversion examples.

\begin{figure}[htb!] 
\begin{center}
\includegraphics[scale=0.5,trim=0cm 0cm 0cm 0cm, clip]{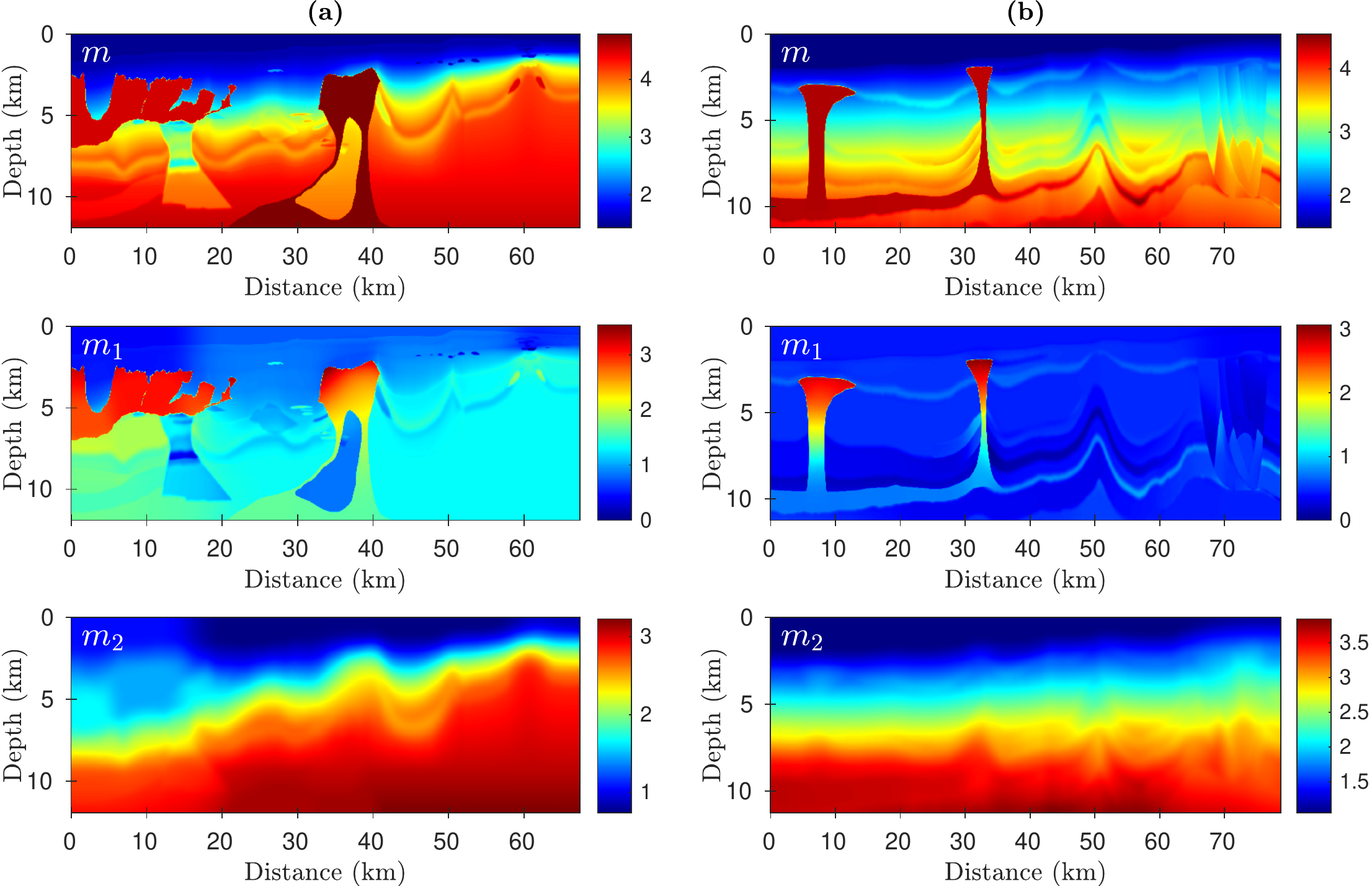}
\caption{Model decomposition into smooth and nonsmooth components. Column (a): the 2004 BP velocity model. Column (b): the 2007 BP velocity model. Top row: the original models $\bm$. Middle row: the nonsmooth component $\bm_1$. Bottom row: the smooth component $\bm_2$.} \label{fig:benchvels}
\end{center}
\end{figure}

\begin{figure}[htb!] 
\begin{center}
\includegraphics[scale=0.6,trim=0cm 0cm 0cm 0cm, clip]{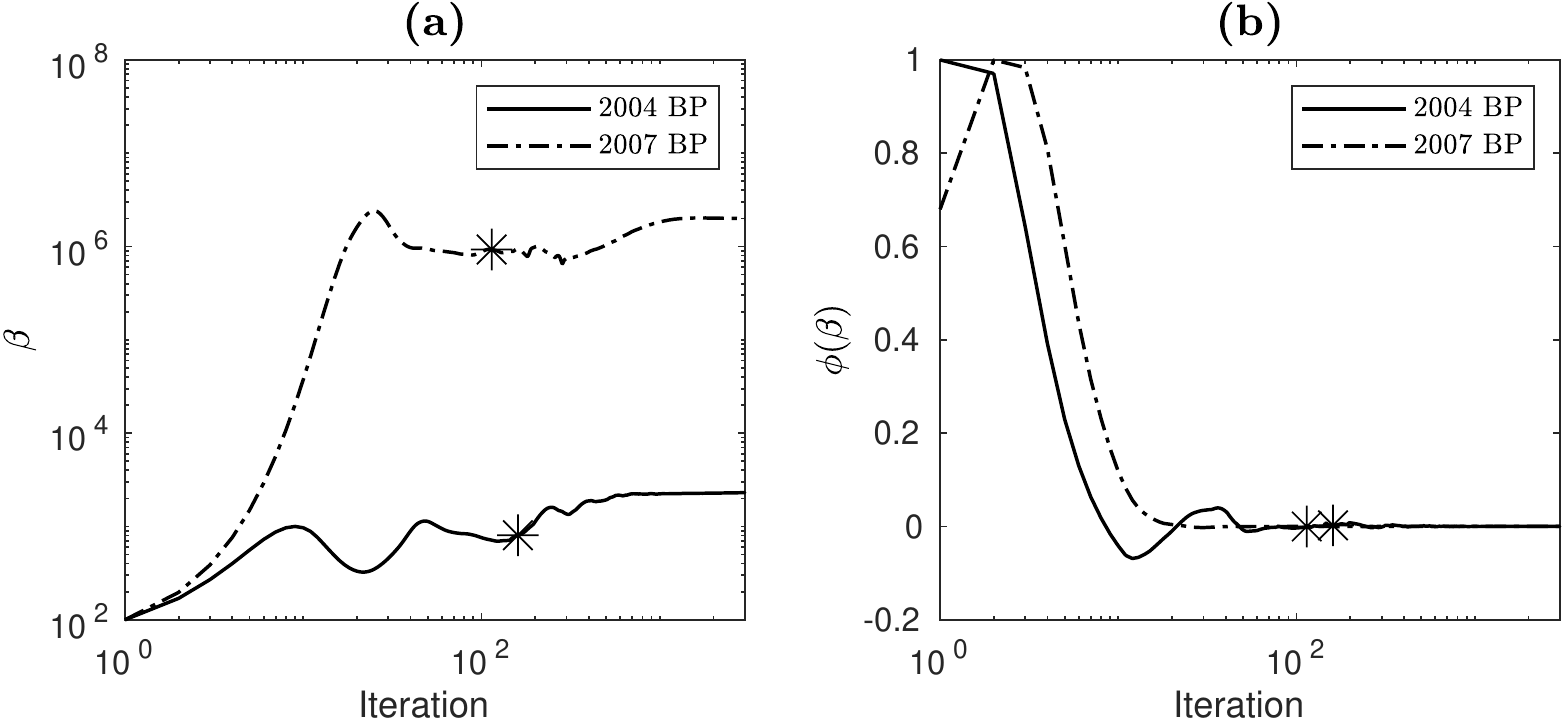}
\caption{Model decomposition into smooth and nonsmooth components. Left frame: values of the balancing parameter $\beta$ versus iteration number. Right frame: values of $\phi(\beta)$ versus iteration number. 
{The asterisk highlights iteration for which the stopping criterion on the stabilization of the solution is satisfied.}
}
\label{fig:benchvels_params}
\end{center}
\end{figure}

Finally, still using two popular datasets in the seismic community, we assess the performance of our algorithm when used for image decomposition. 
Fig.~\ref{fig:benchvels}(a)-(b), top row, show the 2004 BP velocity model \cite{Billette_2004_BPB} and the 2007 BP velocity model \cite{Shah_2007_BAV}, respectively; the velocity variation is a piecewise smooth function of space for both these models. 
We apply 500 iterations of the Tikhonov-TV regularization method to decompose each velocity model into smooth and nonsmooth components: the results are still displayed in Fig. \ref{fig:benchvels}, on the middle and bottom rows, respectively. 
We can clearly observe that the two components $\bold{m}_1$ and $\bold{m}_2$ are optimally separated and allow to recover complementary parts of the original velocity model. 
Fig.~\ref{fig:benchvels_params} shows the evolution of the values of $\beta$ and $\phi(\beta)$ versus iteration for these two test problems. 
We can observe that, in both cases, the simple update scheme in \eqref{fixedPit} stably converges to a root of $\phi(\beta)$, leading to an optimal separation between the smooth and nonsmooth components of the model $\bm$. 

\subsection{X-ray tomography}

\begin{figure}[htb!] 
\begin{center}
\includegraphics[scale=0.7,trim=0cm 0cm 0cm 0cm, clip]{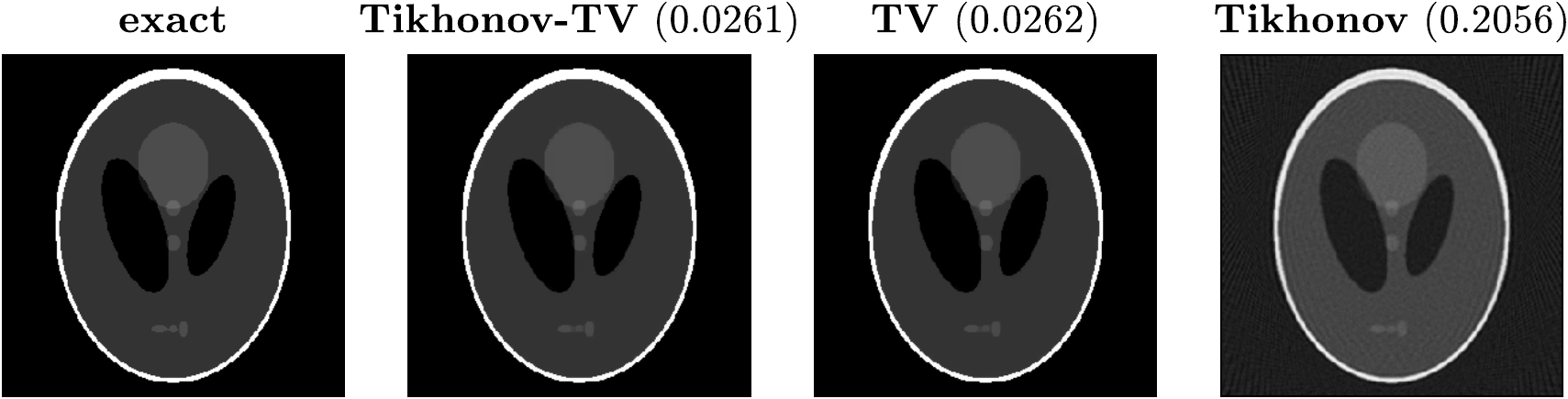}
\caption{X-ray tomography test problem. From left to right: true image (320 $\times$ 320 pixels), estimation by the proposed balanced Tikhonov-TV method, estimation by the TV method, and estimation by the Tikhonov method. The numbers in parenthesis report the relative error of each estimate.
}
\label{fig:phantoms}
\end{center}
\end{figure}

\begin{figure}[htb!] 
\begin{center}
\includegraphics[scale=0.7,trim=0cm 0cm 0cm 0cm, clip]{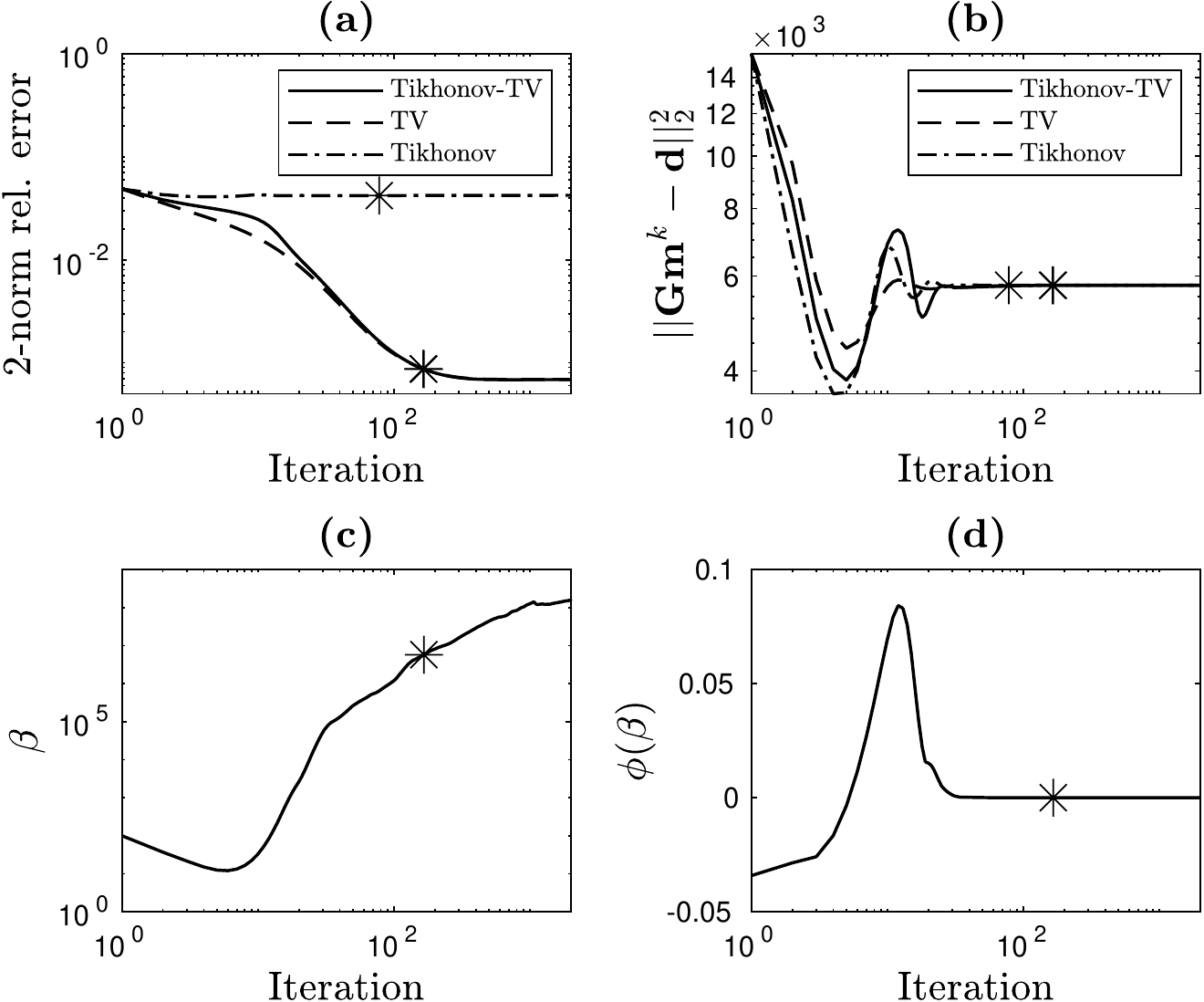}
\caption{X-ray tomography test problem. Evolution of the: (a)  relative errors, (b) squared discrepancy versus iteration for the methods considered in Fig. \ref{fig:phantoms}. Evolution of: (c) $\beta$, (d) $\phi(\beta)$ for the Tikhonov-TV method.
{The asterisk highlights iteration for which the stopping criterion on the stabilization of the solution is satisfied.}
}
\label{fig:phantoms_conv}
\end{center}
\end{figure}

\begin{figure}[htb!] 
\begin{center}
\includegraphics[scale=0.7,trim=0cm 0cm 0cm 0cm, clip]{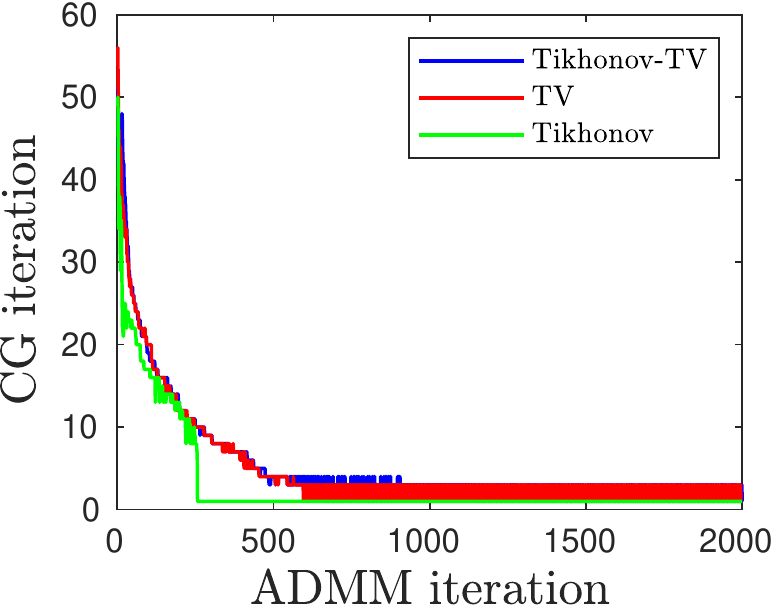}
\caption{X-ray tomography test problem. The number of CG iterations performed at each ADMM iteration.
}
\label{fig:phantoms_CGiter}
\end{center}
\end{figure}

We first consider a parallel X-ray computed tomography test problem. The unknown $\bm$ to be reconstructed is the Shepp-Logan phantom of size $320\times 320$ pixels. 
The original phantom is shown on the leftmost frame of Fig. \ref{fig:phantoms}. We consider a full projection setting with 453 parallel rays for each of the 90 equispaced angles in range $[-90^\circ,+90^\circ]$. This results in a discrete forward operator $\bf G$ of size $40770 \times 102400$. Gaussian white noise of level $1\%$ is added to the data. During the inversion stage, at the $k$th ADMM iteration, the subproblem \eqref{m_sub} for $\bm$ is solved by using the CG method. The CG iterations are stopped when the relative residual tolerance $10^{-7}$ is reached, or when the maximum number of 100 iterations is reached; the estimated solution at the $(k-1)$th ADMM iteration is used as an initial guess for the CG at current iteration. This resulted in a dynamic number of inner CG iterations, as shown in Fig. \ref{fig:phantoms_CGiter}. 
The behavior of other relevant quantities is displayed in {\color{blue}Fig. \ref{fig:phantoms_conv}}. In particular, we can observe that for this test problem, whose solution solely contains piecewise constant features, there are almost no differences in the performance of Tikhonov-TV and TV regularization methods. 

\begin{figure}[htb!] 
\begin{center}
\begin{tabular}{cccc}
\hspace{-1.1cm}\small{\textbf{exact}} & 
\hspace{-1.2cm}\small{\textbf{Tikhonov-TV} (0.2470)} & 
\hspace{-1.2cm}\small{\textbf{TV} (0.2656)} & 
\hspace{-1.2cm}\small{\textbf{Tikhonov} (0.2620)}\\
\hspace{-1.1cm}\includegraphics[scale=0.2]{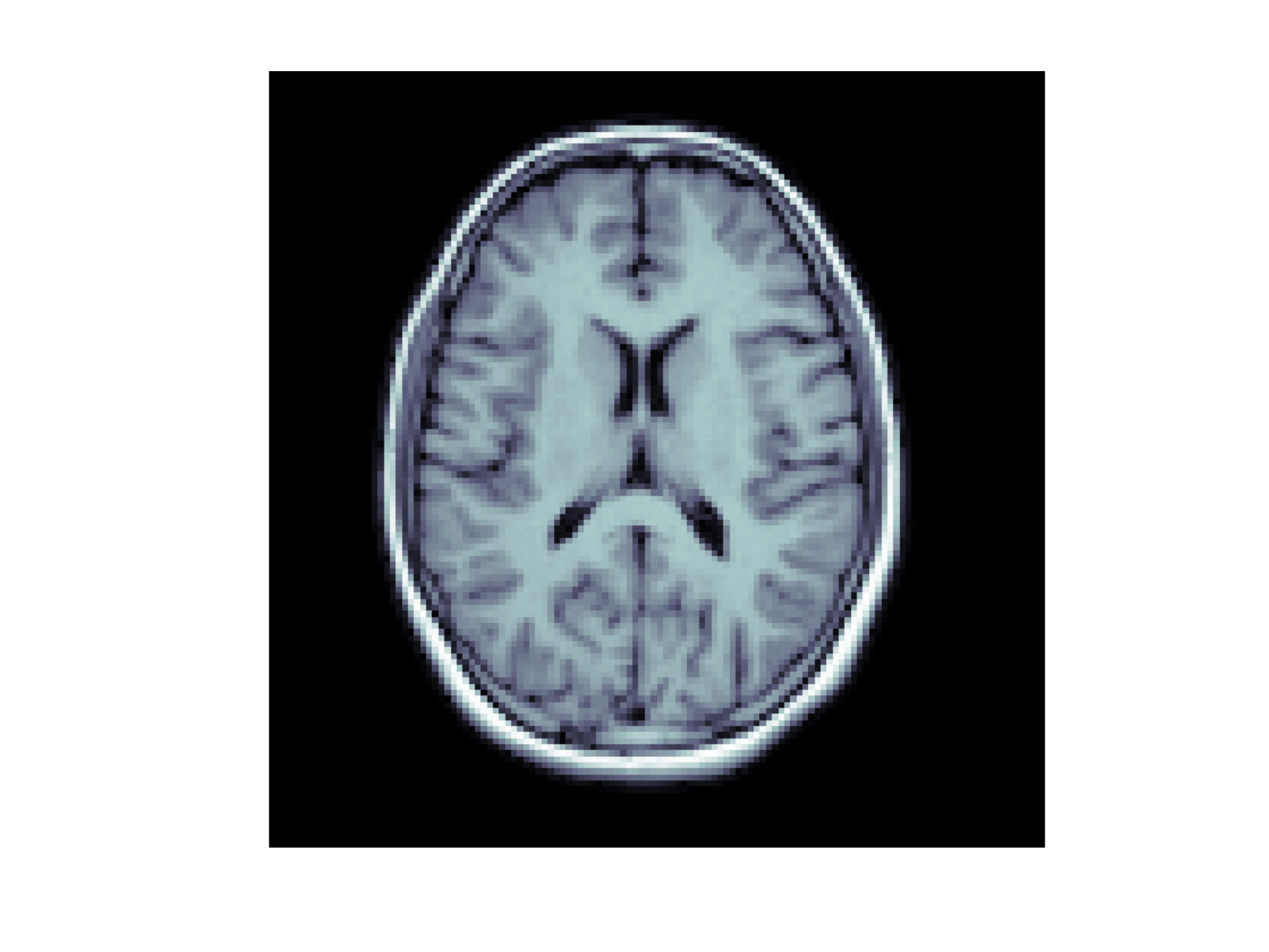} &
\hspace{-1.2cm}\includegraphics[scale=0.2]{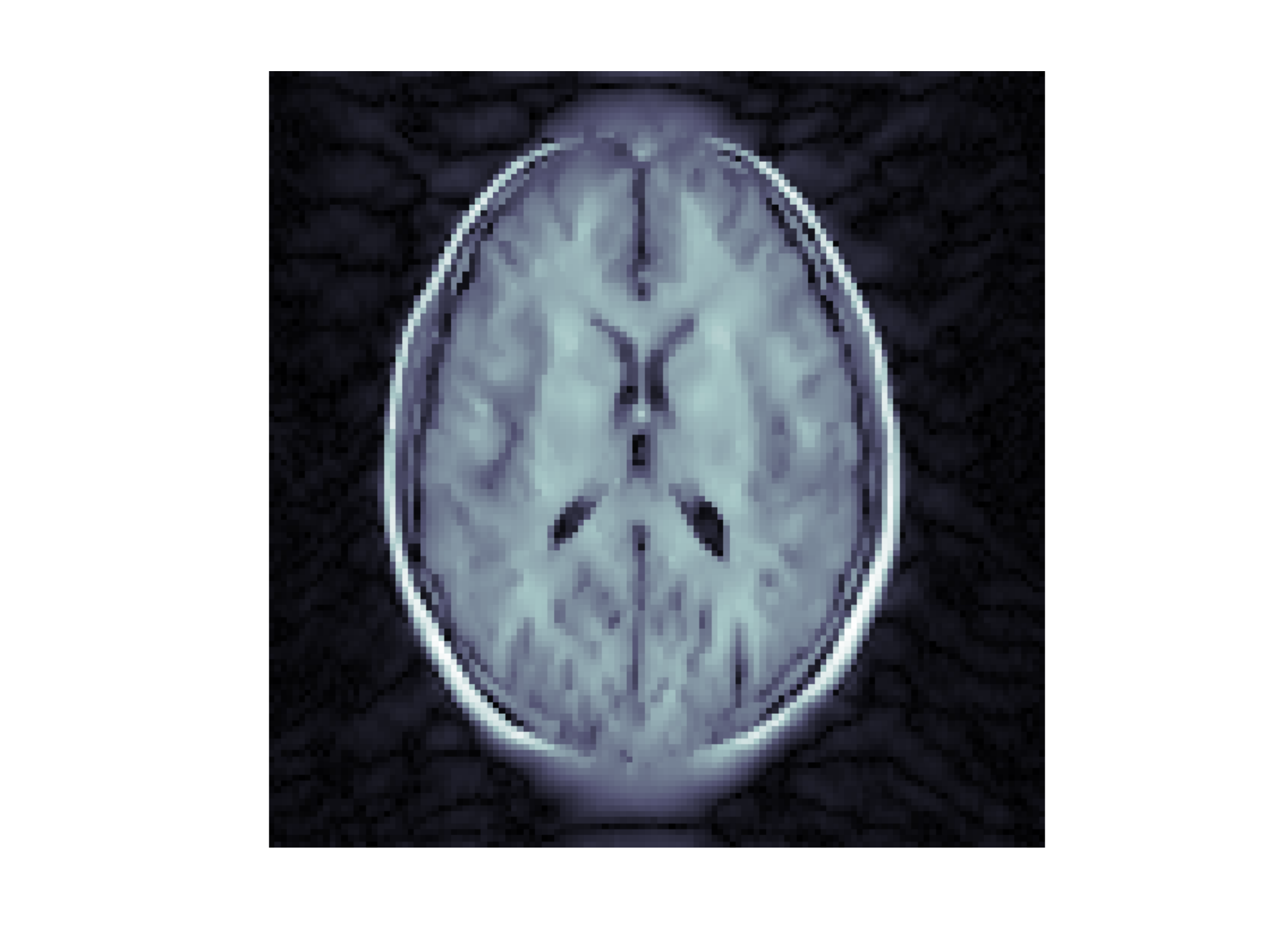} &
\hspace{-1.2cm}\includegraphics[scale=0.2]{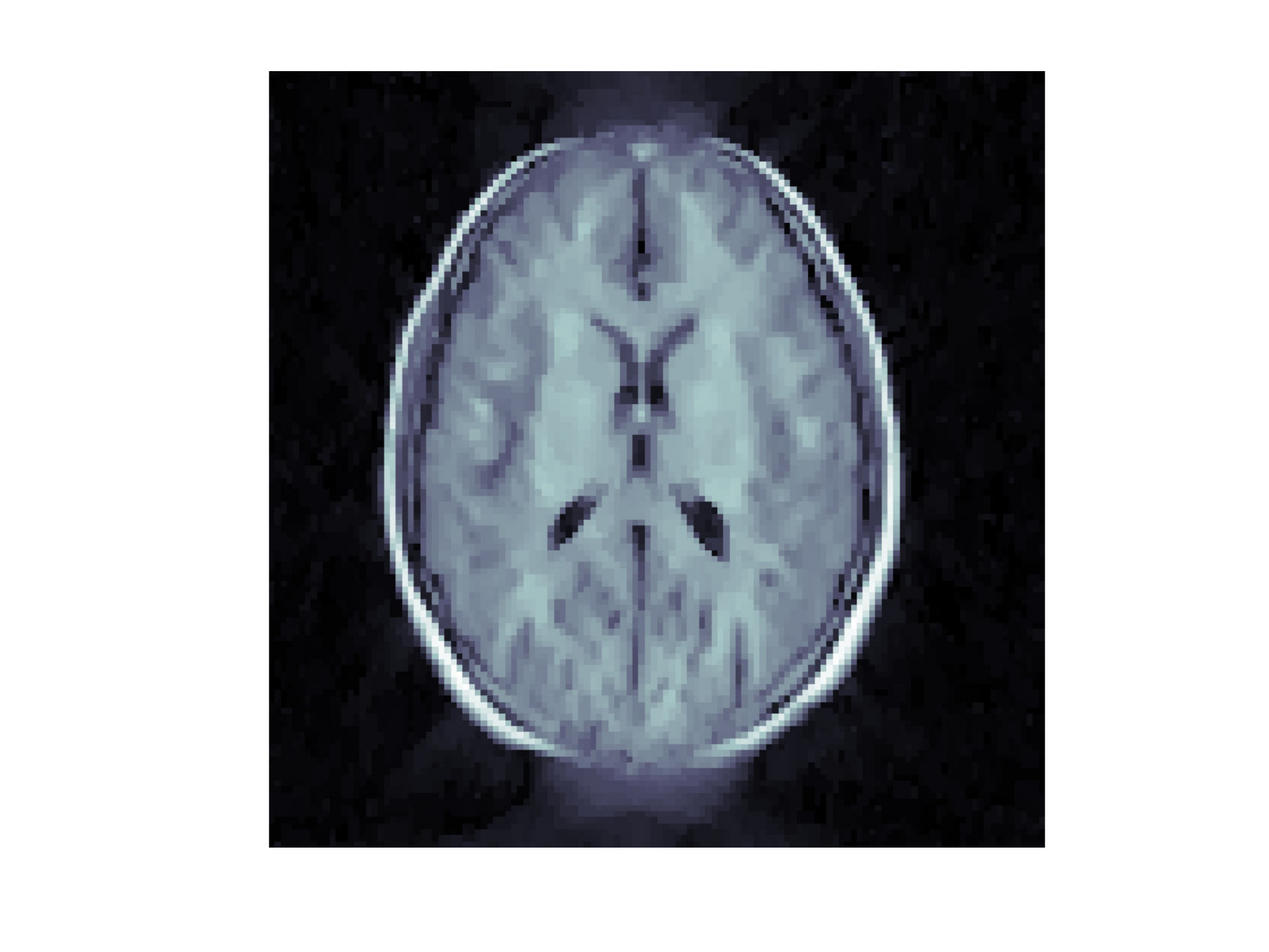} &
\hspace{-1.2cm}\includegraphics[scale=0.2]{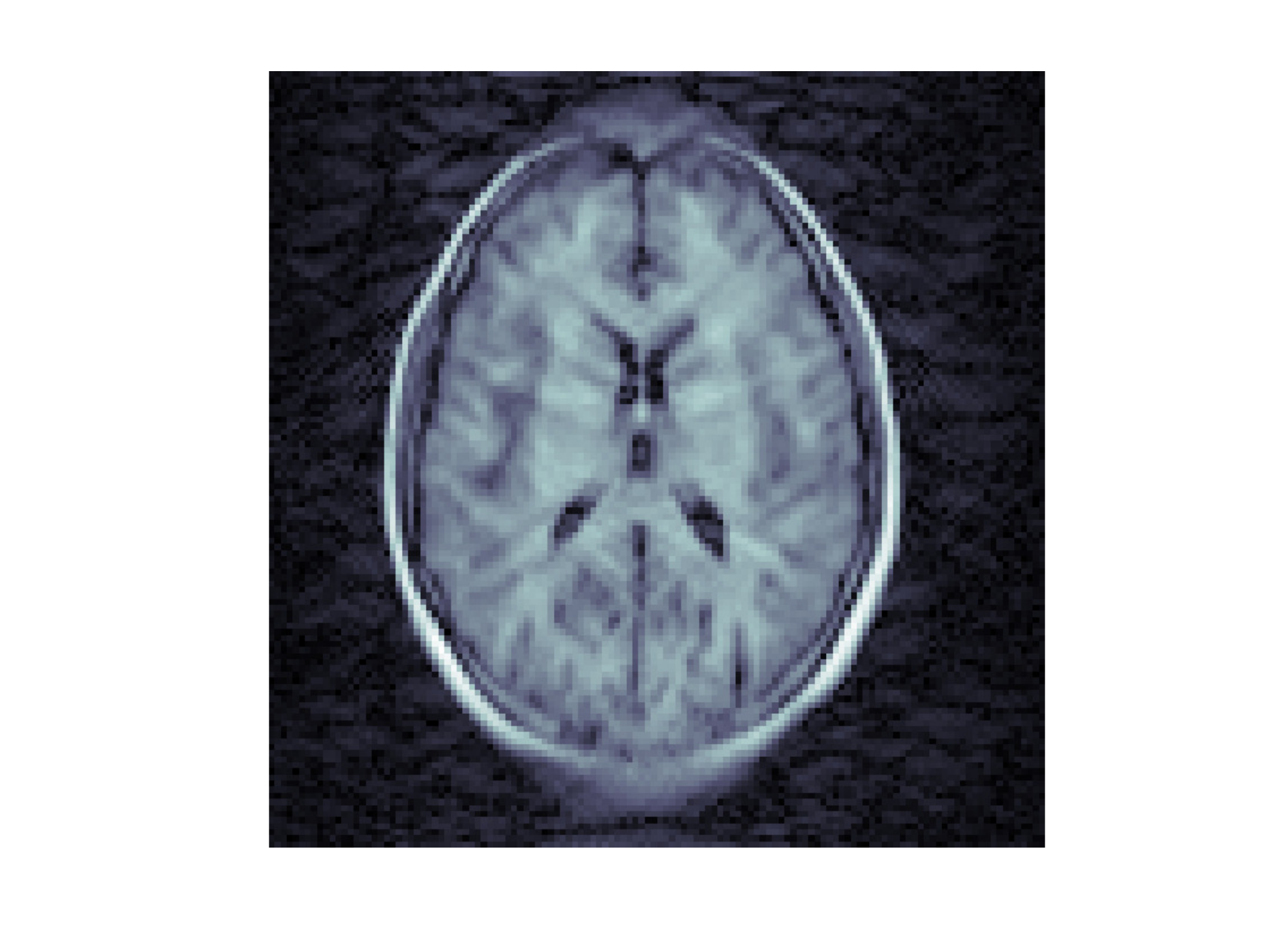}
\end{tabular} 
\caption{X-ray tomography test problem, with limited angles. Exact phantom along with the reconstructions achieved after 600 iterations of each methods. The numbers between brackets are the 2-norm relative errors associated to each reconstruction. To better highlight differences, the square root of the modulus of each pixel is displayed.}
\label{fig:limang_phantoms}
\end{center}
\end{figure}

\begin{figure}[htb!] 
\begin{center}
\includegraphics[scale=0.7,trim=0cm 0cm 0cm 0cm, clip]{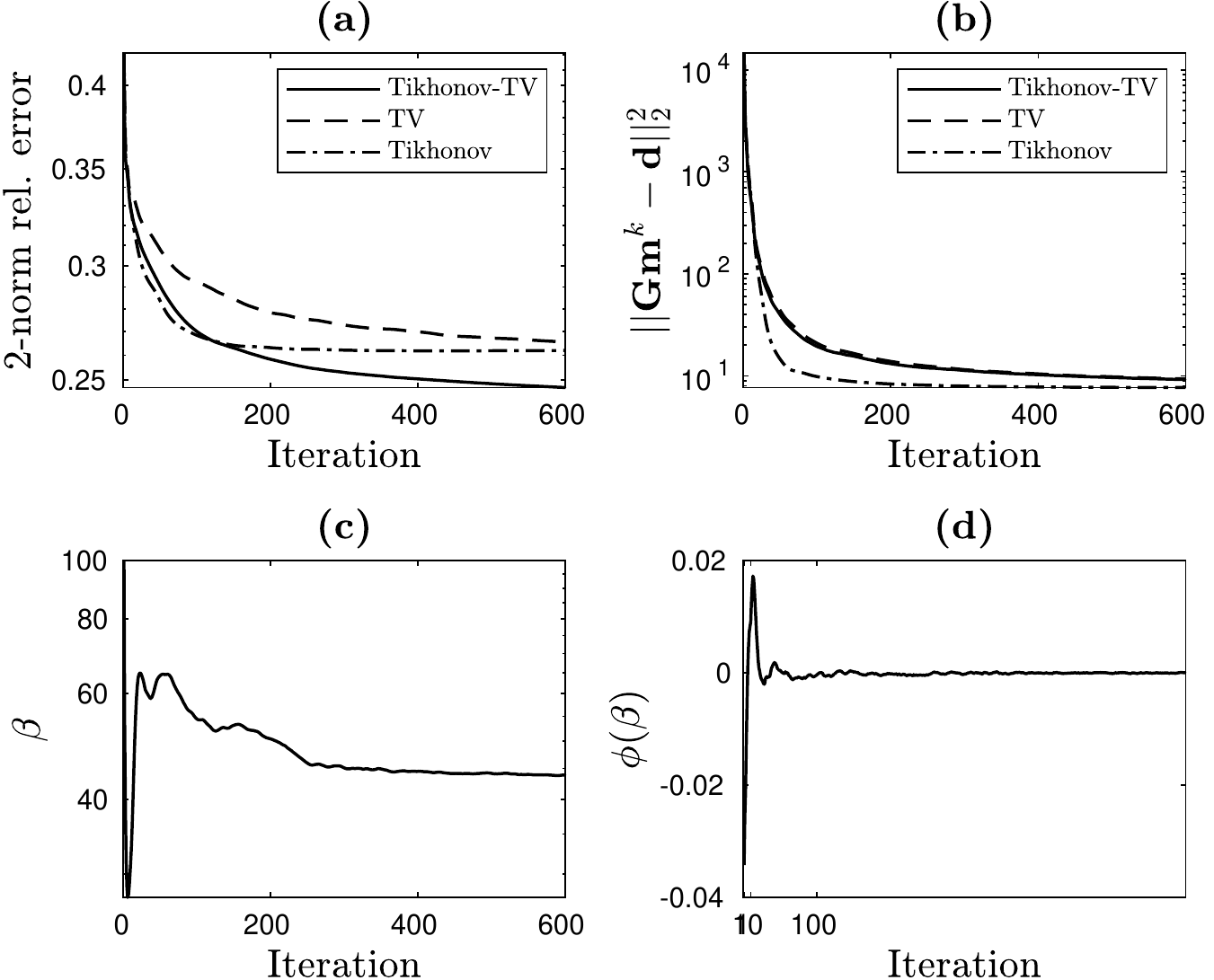}
\caption{X-ray tomography test problem, with limited angles. Frame (a): 2-norm relative errors versus iteration number for different methods. Frame (b): squared 2-norm discrepancies versus iteration number for different methods. Frame (c): values of $\beta$ versus iteration number for Tikhonov-TV. Frame (d): values of $\phi(\beta)$ versus iteration number for Tikhonov-TV.}
\label{fig:limang_graphs}
\end{center}
\end{figure}

We then take a phantom containing both piecewise constant and smooth features, and we resize it to be $128\times 128$ pixels; the resulting phantom is shown on the leftmost frame of Fig.~\ref{fig:limang_phantoms}. We consider a setting where projections along parallel X-rays can be performed only along limited angles: specifically, we consider 181 parallel rays for each of the 85 equispaced angles in $[-42^{\circ}, +42^{\circ}]$. This results in a discrete forward operator $\bold{G}$ of size $15385\times 16384$. Gaussian white noise of level $0.1\%$ is added to the data. For this particular test problem, ADMM applied to solve the Tikhonov-TV and TV-only test problem display a very slow convergence (this is visible in Fig. \ref{fig:limang_graphs} (b), where the only method whose squared residual approximately equals $\varepsilon$ within 600 iterations is Tikhonov regularization). Despite this, the reconstruction computed by the Tikhonov-TV method achieves the lowest relative error among the considered methods, with a good resolution of the background, the piecewise constant and the smooth features: this is clearly visible in Fig.~\ref{fig:limang_phantoms}. {\color{blue}Fig. \ref{fig:limang_graphs}} reports the progress of other relevant quantities versus the iteration count. 

\section{Conclusions} \label{CONS}
In this paper, we proposed a method for selecting the optimal balancing parameter in Tikhonov-TV regularization for the solution of discrete ill-posed inverse problems with piecewise smooth solutions. Namely, the solution of the inverse problem is split into the sum of a smooth and a piecewise constant component, which are separately regularized by a Tikhonov and a TV term that must be balanced. 
We have used robust statistical methods for determining the optimal balance of these terms, motivated by the fact that the gradient entries of the piecewise constant component at jump locations can be considered as anomalies/outliers in comparison to the other gradient entries constituting the smooth background. This led to the characterization of the best balancing parameter as a root of a scalar function. 
Finally, an extremely simple update scheme was proposed for determining the balancing parameter, which can be naturally coupled with {ADMM} to solve the Tikhonov-TV regularized problem. 
Extensive numerical experiments on different inverse problems demonstrate that high quality reconstructions can be obtained by applying the proposed algorithm, and validate the robustness of the proposed selection strategy of balancing parameter.

\begin{appendices}
\section{Proof that the function $\psi(\beta)$ in \eqref{fixedPit} is not expansive.}\label{append}
Let $0<\beta_1<\beta_2$. Let us assume that $\|\nrm(\bg(\beta_1))\|_{\infty}\leq\|\nrm(\bg(\beta_2))\|_{\infty}$. This is a meaningful assumption since $\nrm(\bg(\beta_1))\subset\nrm(\bg(\beta_2))$ because $\bg(\beta_1)$ is expected to be less smooth than $\bg(\beta_2)$ (because the smooth component $\bg_2$ of $\bg$ is weighted less in the objective function \eqref{main_con2}). We know that $\|\bg_2(\beta_1)\|_2>\|\bg_2(\beta_2)\|_2$ (still as a consequence of the weighting in the objective function \eqref{main_con2}). Since the 2-norm and $\infty$-norm are equivalent on $\mathbb{R}^N$ (namely, for all $\bold{x}\in\mathbb{R}^N$, $\|\bold{x}\|_{\infty}\leq \|\bold{x}\|_2\leq N^{1/2}\|\bold{x}\|_{\infty}$), the 2-norm inequality above implies the $\infty$-norm inequality $\|\bg_2(\beta_1)\|_{\infty}\geq N^{-1/2}\|\bg_2(\beta_2)\|_{\infty}$. To keep the notations light, in the following we let
\[
n_1:=\|\nrm(\bg(\beta_1))\|_{\infty},\; n_2:=\|\nrm(\bg(\beta_2))\|_{\infty},\; d_1:=\|\bg_2(\beta_1)\|_{\infty},\;d_2:=\|\bg_2(\beta_2)\|_{\infty}.
\]
As a consequence of the inequalities above, which can be rewritten as
\begin{equation}\label{ineqfixP}
n_1\leq n_2\quad\mbox{and}\quad d_1\geq N^{-1/2}d_2\quad\left(\mbox{so that}\quad \frac{n_1}{d_1}<\frac{n_2}{N^{-1/2}d_2}\right)\,,
\end{equation}
the term appearing in the right-hand side parenthesis in \eqref{fixedPit} can be 
rewritten and bounded as follows
\begin{eqnarray*}
\left(4+\frac{n_1}{d_1}\right)^{-1}-1&=&N^{-1/2}\left(\left(4N^{-1/2}+\frac{n_1}{N^{1/2}d_1}\right)^{-1}-N^{1/2}\right)\\
&>&N^{-1/2}\left(\left(4N^{-1/2}+\frac{n_2}{d_2}\right)^{-1}-N^{1/2}\right)\\
&\geq&N^{-1/2}\left(\left(4+\frac{n_2}{d_2}\right)^{-1}-N^{1/2}\right)\,,
\end{eqnarray*}
where, in the first inequality, we have exploited \eqref{ineqfixP}. Therefore, the modulus of the difference $\psi(\beta_1) - \psi(\beta_2)$ can be rewritten and bounded as
\begin{eqnarray*}
\left|\psi(\beta_1)-\psi(\beta_2)\right| &=&\left| \left(\left(4+\frac{n_1}{d_1}\right)^{-1}-1\right)\beta_1-\left(\left(4+\frac{n_2}{d_2}\right)^{-1}-1\right)\beta_2\right|\\
&<& \left| \left(\left(4+\frac{n_1}{d_1}\right)^{-1}-1\right)\beta_1
-\left(N^{-1/2}\left(4+\frac{n_2}{d_2}\right)^{-1}-1\right)\beta_2\right| \\
&\leq& \left| \left(\left(4+\frac{n_1}{d_1}\right)^{-1}-1\right)\beta_2-\left(\left(4+\frac{n_1}{d_1}\right)^{-1}-1\right)\beta_1\right|\\
&=& \left| \left(\left(4+\frac{n_1}{d_1}\right)^{-1}-1\right)\right| |\beta_1-\beta_2|\,.
\end{eqnarray*}
Since the first factor on the left-hand side is positive and strictly less than 1, we can conclude that $\psi(\beta)$ is not expansive. 
\end{appendices}

\end{document}